\documentclass[a4paper,11pt,reqno]{article}
\usepackage{amsmath,amsfonts,amssymb,amsthm}
\usepackage{mathrsfs}
\usepackage{color}
\usepackage{graphicx}
\usepackage[numbers,sort&compress]{natbib}
\usepackage[colorlinks=true]{hyperref}
\hypersetup{urlcolor=blue,citecolor=blue}

\def\R{{\mathbb{R}}}
\numberwithin{equation}{section}

\begin{document}

\title{On the two-component $b$-family system modeling shallow-water waves of medium amplitude}

\author{Mengwen $\mbox{Huang}$\footnote{E-mail:
ts26080008a31tm@cumt.edu.cn}\quad
and \quad Xingxing $\mbox{Liu}$\footnote{Corresponding author. E-mail: liuxxmaths@cumt.edu.cn}\\
$\mbox{School}$ of mathematics, China University of
Mining and Technology,\\Xuzhou, Jiangsu 221116, China}

\date{}
\maketitle

\begin{abstract}

We formally derive a two-component $b$-family system from the full governing equations for two-dimensional
water waves in the medium-amplitude regime. The system contains higher-order nonlinearities and includes
a two-component Camassa-Holm type system for equatorial shallow-water waves with the Coriolis effect as
a special case. We establish the local well-posedness of the Cauchy problem in
$B^{s-1}_{p,r}\times B^s_{p,r}$ with $s>\max\{1+\frac1p,\frac32\}$. Using Littlewood-Paley theory and
Osgood's lemma, we also prove local well-posedness in the critical Besov space
$B^{\frac12}_{2,1}\times B^{\frac32}_{2,1}$. Furthermore, we derive a blow-up criterion for strong solutions.
For a special choice of the parameters, we prove finite-time blow-up for a class of initial data and provide
an upper bound for the maximal existence time.


\end{abstract}

\section{Introduction}
\newtheorem{theorem1}{Theorem}[section]
\newtheorem{lemma1}{Lemma}[section]
\newtheorem {remark1}{Remark}[section]
\newtheorem {definition1}{Definition}[section]
\newtheorem{corollary1}{Corollary}[section]
\newtheorem{proposition1}{Proposition}[section]
\par

Water waves are familiar natural phenomena that also give rise to rich and diverse mathematical structures, making their study perpetually fascinating.
However, the governing equations for water waves, which feature a free boundary and nonlinear convection,
are so challenging that gaining a clear understanding of wave dynamics directly from these equations is often difficult. Consequently, the derivation of
simplified asymptotic model equations has become an important direction in the theory of water waves.

In the non-dimensionalisation of the water-wave problem, two key dimensionless parameters naturally arise: the amplitude parameter
$\varepsilon=a/h_{0}$ (where $a$ is a typical wave amplitude and $h_0$ the mean water depth) and the shallowness parameter $\mu=h_{0}^{2}/\lambda^{2}$ (with $\lambda$ a
characteristic wavelength). It is important to emphasize that the two parameters $\varepsilon$ and $\mu$
are physically independent, and there is no a priori functional relation between them.
Hence, statements such as "$\varepsilon$ is smaller than $\mu$" or "$\varepsilon=O(\mu)$" are not intrinsically meaningful,
as they implicitly assume some relationship between the two parameters.
Nevertheless, the apparent "balance" between $\varepsilon$ and $\mu$ has excited a lot of interest \cite{CJ1}.

Indeed, the relative size of these parameters distinguishes different physical regimes.
In the long-wave regime, known as the Boussinesq scaling ($\mu\ll1$ and $\varepsilon=O(\mu)$),
nonlinearity and dispersion are of the same order, giving rise to the classical Korteweg-de Vries (KdV) and Benjamin-Bona-Mahony (BBM) equations \cite{BBM,KdV}.
However, this scaling requires small amplitudes and cannot capture strong nonlinear phenomena such as wave breaking.
In the regime of medium-amplitude waves, also called the Camassa-Holm scaling ($\mu\ll1$ and $\varepsilon=O(\sqrt{\mu})$),
the Camassa-Holm (CH) \cite{CH} and Degasperis-Procesi (DP) equations \cite{DP} have been derived as valid asymptotic approximations to the full governing equations for water waves \cite{CL,J1,J2,L}.
These two equations exhibit rich dynamics, including peaked solitary waves (peakons) \cite{CH,CS1,LL} and finite-time wave breaking \cite{CE,ELY,LY}, and have therefore attracted considerable interest.
The two classical regimes mentioned above do not exhaust all possible regimes.
As research into strongly nonlinear waves has progressed, it has become clear that higher-order nonlinear terms are necessary to describe their evolution accurately. In recent years, several new asymptotic scalings have been proposed for larger amplitudes. For instance, taking $\varepsilon=O(\sqrt[5]{{\mu}^2})$ or $\varepsilon=O(\sqrt[3]{\mu})$
(where $\varepsilon$ is larger than in the CH scaling) leads to shallow-water models that incorporate cubic or even higher-order nonlinearities \cite{CHL,LLL}.

The success of the aforementioned single-component equations (KdV, CH, etc.) in shallow-water waves studies has naturally inspired the search for generalizations with richer properties.
In particular, two-component shallow-water wave systems have attracted considerable attention because they couple the evolution of the surface elevation and the horizontal velocity.
Consider the classical two-component Camassa-Holm (2CH) system, which has been derived in the theory of water waves under the Boussinesq scaling
on several occasions. For instance, Constantin and Ivanov \cite{CI} established this system using the approach based on Green--Naghdi approximation in the shallow water framework without vorticity.
Subsequently, Ivanov \cite{I} derived a constant-vorticity analogue of this system from the Euler equations for water waves with an underlying constant vorticity. This analogue contains an additional linear term due to the vorticity and reduces to the system of Constantin and Ivanov in the irrotational case.
Moreover, Escher et al. \cite{EHKL} also started from the Euler equations and derived a two-component $b$-family system with constant vorticity, which includes the 2CH system (for $b=2$)
and the two-component Degasperis-Procesi (2DP) system (for $b=3$) as special cases.
Although these derivations follow different paths (one based on the Green--Naghdi approximation \cite{CI} and the others starting directly from the Euler equations \cite{EHKL,I}), they are all
carried out in the small-amplitude regime ($\mu\ll1$ and $\varepsilon=O(\mu)$).

One of the goals of the present paper is to adopt the medium-amplitude regime ($\mu\ll1$ and $\varepsilon=O(\sqrt{\mu})$) and derive formally from the full governing equations for two-dimensional water waves
the following more general two-component $b$-family system.
\begin{equation}\label{1.1}
\left\{
 \begin{aligned}
 &\rho_t+u\rho_x+(b-1)\rho u_{x}-\gamma_1(u^3)_x\rho=0,\\
&m_t+\alpha(um_x+b u_{x}m)+\frac{1}{2}(\rho^2)_x+\gamma_2(u^2\rho^2)_x+\gamma_3(u^2\rho)_x+\gamma_4(u^3)_x=0,\quad m=u-u_{xx},
\end{aligned}
\right.
\end{equation}
where $u(t,x)$ denotes the average horizontal fluid velocity, $\rho(t,x)$ is related to the free
surface elevation from equilibrium, $b>1$ is a real parameter, and the constants
\begin{equation*}\label{1.2}
\gamma_1=\frac{(2k_1k_2-2k_2)(b-1)}{k_1^3}=(2b-3)^2(b-1),\quad
\gamma_2=\frac{\tilde\beta}{k_1^2},\quad
\gamma_3=-\frac{\tilde\beta}{k_1^2\sqrt{b-1}},\quad
\gamma_4=\frac{\beta-\tilde\beta}{k_1^2},
\end{equation*}
with
\begin{equation*}\label{1.3}
k_1=\frac{1}{2(b-1)},\quad
k_2=-\frac{2b-3}{8(b-1)^2},\quad
k_3=\frac{1}{12(b-1)},\quad
\beta=\frac{(2b-3)(5b-7)}{12(b-1)^2},
\end{equation*}
and $\tilde\beta\in\R$ a free parameter, while the coefficient of the transport part is
\begin{equation*}\label{1.3a}
\alpha=\frac{k_1-2k_2}{k_1^2(1+b)}=\frac{4b-5}{b+1},
\end{equation*}
in particular $\alpha=1$ for $b=2$. Throughout this paper we assume that $\alpha\neq 0$
(i.e., $b\neq\frac54$).

Our model system \eqref{1.1} contains several free parameters and admits interesting special cases. For $b=2$ and $\tilde\beta=-\frac34$, one has $\alpha=1$, $\gamma_1=1$, $\gamma_2=-3$,
$\gamma_3=3$, $\gamma_4=4$, and \eqref{1.1}
reduces exactly to
the 2CH type system under the CH scaling \cite{HL,ZL}
 \begin{equation}\label{2CH}
\left\{
 \begin{aligned}
&{\rho _t}  + {(\rho u)_x} -\rho {({u^3})_x} =  0,\\
&{m_t} +  (2m{u_x} + u{m_x})+ \rho\rho_x-3{({u^2}\rho (\rho  - 1))_x} + 4{({u^3})_x} = 0,\ m=u-u_{xx},
\end{aligned}
\right.
\end{equation}
which models equatorial shallow-water waves of medium amplitude with the Coriolis effect
\cite{HL}. Compared with derivations based on the Green--Naghdi approximation \cite{HL,ZL}, our approach relies directly on the Euler equations and retains higher-order non-local nonlinear terms in the medium-amplitude regime. Moreover, it provides a more direct framework for subsequently incorporating complex factors such as rotation and vorticity.

The other goal of this paper is to investigate the Cauchy problem for (\ref{1.1}). To handle the higher-order
nonlinearities and the coupling between the two components, we adopt a unified transport-nonlocal formulation
and use Littlewood-Paley theory. Theorem \ref{LWP} establishes local existence and uniqueness in
$B^{s-1}_{p,r}\times B^s_{p,r}$ with $s>\max\{1+\frac1p,\frac32\}$, together with H\"older continuous
dependence on the initial data in suitable lower-regularity spaces. At the critical regularity, where the
usual contraction argument does not close directly, we further combine a logarithmic interpolation inequality
with Osgood's lemma to establish local well-posedness in $B^{\frac12}_{2,1}\times B^{\frac32}_{2,1}$ in
Theorem \ref{LWPcri}. In particular, for $b=2$ and $\tilde\beta=-\frac34$, Theorem \ref{LWP} recovers the
local well-posedness result of Hu and Liu \cite{HL} for (\ref{2CH}) in $H^{s-1}\times H^s$ with
$s>\frac32$, while Theorem \ref{LWPcri} also provides a critical-space local theory for this special system.

For the blow-up analysis, we consider two complementary questions: how to characterize continuation of a
general strong solution, and how to exhibit initial data that actually lead to finite-time blow-up. Theorem
\ref{blowup} answers the first question through a necessary and sufficient condition involving the time
integral of certain $L^\infty$-norms of $u$, $u_x$ and $\rho$. The criterion follows directly from a priori
estimates and does not require a special conservation law. Since such a criterion does not itself produce a
blow-up example, we further consider the parameters $b=3$ and $\tilde\beta=\frac12$. In this case, the
preservation of odd $u$ and even $\rho$ reduces the blow-up argument at $x=0$ to a
Riccati-type differential inequality for $u_x(t,0)$. Theorem \ref{symmetryblowup} thereby exhibits a concrete
class of initial data for which finite-time blow-up occurs and gives an explicit upper bound for the maximal
existence time.

The rest of the paper is organized as follows. In Section 2, we derive (\ref{1.1}) formally from the
full water-wave equations and discuss its asymptotic simplifications. Section 3 is devoted to the local
well-posedness in both noncritical and critical Besov spaces. In Section 4, we establish the general blow-up
criterion and prove finite-time blow-up for the special class of initial data described above. The Appendix
collects the definitions and auxiliary estimates used in the proofs.

\section{Derivation of the model system}\label{sec2}
\newtheorem{theorem2}{Theorem}[section]
\newtheorem{lemma2}{Lemma}[section]
\newtheorem {remark2}{Remark}[section]
\newtheorem {definition2}{Definition}[section]
\newtheorem{corollary2}{Corollary}[section]
\par
In this section, we work from the dimensionless governing equations for two-dimensional water waves and carry out an asymptotic analysis in the medium-amplitude regime
($\mu\ll1$ and $\varepsilon=O(\sqrt{\mu})$).
Unlike previous derivations based on the Green--Naghdi approximation \cite{CI,HL,ZL}, our approach follows the Euler equations route of \cite{EHKL,I} and relies directly on the governing equations for water waves.
By systematically expanding the velocity field and pressure, and introducing suitable auxiliary variables, we derive a two-component $b$-family shallow-water system valid for medium-amplitude waves, namely the system (\ref{system}) below.

As a starting point, we present the dimensionless Euler equations together with the corresponding boundary conditions
\begin{align}\label{cases1}
\begin{cases}
u_{t}+\varepsilon(uu_{x}+wu_{z})=-p_{x}, &\text{in } 0<z<1+\varepsilon\eta(t,x),\\
\mu\big(w_{t}+\varepsilon(uw_{x}+ww_{z})\big)=-p_{z}, &\text{in } 0<z<1+\varepsilon\eta(t,x),\\
u_{x}+w_{z}=0, &\text{in } 0<z<1+\varepsilon\eta(t,x),\\
u_{z}-\mu w_{x}=0, &\text{in } 0<z<1+\varepsilon\eta(t,x),\\
p=\eta, &\text{on } z=1+\varepsilon\eta(t,x),\\
w=\eta_{t}+\varepsilon u\eta_{x}, &\text{on } z=1+\varepsilon\eta(t,x),\\
w=0, &\text{on } z=0.\\
\end{cases}
\end{align}
Here $u,w,p,\eta$ denote the scaled and nondimensionalised horizontal and vertical velocity components, pressure, and free-surface elevation, respectively, for which we retain the same symbols as for their dimensional counterparts.
The parameters $\varepsilon=a/h_{0}$ and $\mu=h_{0}^{2}/\lambda^{2}$ are defined as in Introduction.
The derivation of (\ref{cases1})
from the Euler equations via the standard nondimensionalisation and scaling can be found in \cite{CJ1,J1,L}.

From the third, fourth and seventh equations of (\ref{cases1}), we obtain
\begin{equation}\label{u}
u = {u_0} - \mu \frac{{{z^2}}}{2}{u_{0,xx}} + O({\mu ^2}),
\end{equation}
and
\begin{equation}\label{w}
w =  - z{u_{0,x}} + \mu \frac{{{z^3}}}{6}{u_{0,xxx}} + O({\mu ^2}),
\end{equation}
where $u_0(t,x)$ is the leading-order approximation for $u$ in an asymptotic expansion. Note that $u_0$ is independent of $z$, i.e., $u_0 = u_0(t,x)$, since from the fourth equation of (\ref{cases1})
it follows that $u_z=0$ when $\mu\rightarrow 0$.
Substituting $u$ and $w$ into the second equation of (\ref{cases1}), we get
\begin{equation*}
{p_z} = \mu z({u_{0,tx}} + \varepsilon {u_0}{u_{0,xx}} -\varepsilon u_{0,x}^2) +O({\mu ^2}),
\end{equation*}
which along with the boundary condition of $p$ on $z=h=1+\varepsilon\eta$ implies
\[p = \eta  - \frac{\mu }{2}({h^2} - {z^2})({u_{0,tx}} + \varepsilon {u_0}{u_{0,xx}} - \varepsilon u_{0,x}^2) + O({\mu ^2}).\]
Differentiating the above equation with respect to $x$, we have
\begin{equation*}
\begin{split}
{p_x} ={\eta _x} - \frac{\mu }{2}({h^2} - {z^2}){({u_{0,tx}} + \varepsilon {u_0}{u_{0,xx}} - \varepsilon u_{0,x}^2)_x} - \mu h{h_x}({u_{0,tx}}
+ \varepsilon {u_0}{u_{0,xx}} - \varepsilon u_{0,x}^2)+ O({\mu ^2}).
\end{split}
\end{equation*}
On the other hand, in view of the expressions of $u,w$ in (\ref{u})-(\ref{w}), we deduce from the first equation of (\ref{cases1}) that
\begin{equation*}
\begin{split}
{u_{0,t}} - \mu \frac{{{z^2}}}{2}{u_{0,txx}} + \varepsilon {u_0}{u_{0,x}}+ \varepsilon \mu \frac{{{z^2}}}{2}({u_{0,x}}{u_{0,xx}} - {u_0}{u_{0,xxx}}) =-{p_x} + O({\mu ^2}).
\end{split}
\end{equation*}
Combining the last two equations, we find
\begin{equation}\label{u0}
\begin{split}
{u_{0,t}}+ {\eta _x}+ \varepsilon {u_0}{u_{0,x}}= \frac{\mu }{2}\big({h^2}{({u_{0,tx}} + \varepsilon {u_0}{u_{0,xx}} - \varepsilon u_{0,x}^2)}\big)_x+ O({\mu ^2}).
\end{split}
\end{equation}
In the medium-amplitude regime ($\mu\ll1$ and $\varepsilon=O(\sqrt{\mu})$), we can expand (\ref{u0}) up to the order of $O(\varepsilon^2\mu,\mu ^2)$ as follows.
\begin{equation}\label{u0final}
\begin{split}
(u_{0}-\frac{\mu }{2}u_{0,xx})_t + {\eta _x}+ \varepsilon {u_0}{u_{0,x}}= \frac{\varepsilon\mu }{2}( {u_0}{u_{0,xx}} -  u_{0,x}^2+2\eta u_{0,tx})_x
+ O({\varepsilon^2\mu,\mu ^2}).
\end{split}
\end{equation}
Moreover, combining the sixth equation of (\ref{cases1}) with the expressions of $u,w$ in (\ref{u})-(\ref{w}), we obtain
\begin{equation}\label{first1}
- (1+\varepsilon\eta){u_{0,x}} + \mu \frac{{{(1+\varepsilon\eta)^3}}}{6}u_{0,xxx} =\eta_t+\varepsilon \big( {u_0} - \mu \frac{{{(1+\varepsilon\eta)^2}}}{2}{u_{0,xx}} \big)\eta_x+O({\mu ^2}).
\end{equation}
Thus, we derive an equation containing all terms up to the order $ O({\varepsilon^2\mu,\mu ^2})$ from (\ref{first1}) that
\begin{equation}\label{first}
\eta_t+\big((1+\varepsilon\eta)u_0\big)_x-\frac{1}{6}\mu u_{0,xxx}-\frac{1}{2}\varepsilon\mu(\eta u_{0,xx})_x=O({\varepsilon^2\mu,\mu ^2}).
\end{equation}

Up to the correction $O(\varepsilon,\mu),$ it is thus inferred from (\ref{u0final}) and (\ref{first}) that
\begin{equation}\label{linearueta}
\left\{
\begin{aligned}
&{\eta _t}+u_{0,x}=O(\varepsilon,\mu),\\
&{u_{0,t}}+{\eta _x}=O(\varepsilon,\mu),
\end{aligned} \right.
\end{equation}
giving us
\begin{equation*}
\left\{
\begin{aligned}
&{\eta _{tt}}  - {\eta _{xx}} = O(\varepsilon ,\mu ),\\
&{u_{0,tt}} - {u_{0,xx}} = O(\varepsilon ,\mu ).
\end{aligned}
\right.
\end{equation*}
Then we get the general solution of the above system
\begin{equation*}
\left\{
\begin{aligned}
&\eta  = \eta_1 (x - t) + \eta_2 (x + t) + O(\varepsilon ,\mu ),\\
&u_0= u_{0,1}(x - t) + u_{0,2} (x + t) + O(\varepsilon ,\mu ).
\end{aligned}
\right.
\end{equation*}
It is convenient first to restrict ourselves to waves moving towards the right side, i.e.,
\begin{equation*}
\left\{
\begin{aligned}
&\eta  = \eta (x - t) + O(\varepsilon ,\mu ),\\
&u_0 = u_0(x - t) + O(\varepsilon ,\mu ),
\end{aligned}
\right.
\end{equation*}
Thus we have
\begin{equation}\label{linearueta1}
\left\{
 \begin{aligned}
& {\eta _t} =  - {\eta _x} + O(\varepsilon ,\mu ),\\
&{u_{0,t}} =  - {u_{0,x} }+ O(\varepsilon ,\mu ).
\end{aligned}
\right.
\end{equation}
It then follows from (\ref{linearueta}) and (\ref{linearueta1}) that
\begin{equation}\label{relationueta}
{u_{0,t}} + {\eta _x} = O(\varepsilon ,\mu ),\quad\quad       u_0 - \eta  = O(\varepsilon ,\mu ),\quad\quad
{u_{0,x}} +{\eta _t} = O(\varepsilon ,\mu ).
\end{equation}

To derive our model system, we introduce the auxiliary variable $\varrho$, defined in terms of $\eta$ and $u_0$ as follows
\begin{equation}\label{varrho}
\varrho=1+k_1\varepsilon\eta+k_2\varepsilon^2(u_0^2+\eta^2)+k_3\varepsilon\mu u_{0,xx},
\end{equation}
for some constants $k_1,k_2$ and $k_3$. As we proceed, it will be seen that these constants can be determined (see (\ref{relation1})-(\ref{relation2}) and (\ref{relation3}) below).
Combining (\ref{u0final}), (\ref{first}) with (\ref{relationueta})-(\ref{varrho}), we have
\begin{eqnarray}\label{varrhot}
{\varrho _t}&=&k_1\varepsilon\eta_t+2k_2\varepsilon^2(u_0u_{0,t}+\eta\eta_t)+k_3\varepsilon\mu u_{0,txx}\nonumber \\
&=&- k_1\varepsilon\big((1+\varepsilon\eta)u_0\big)_x+\frac{k_1}{6}\varepsilon\mu u_{0,xxx}-2k_2\varepsilon^2 u_0({\eta _x}+ \varepsilon {u_0}{u_{0,x}})\nonumber \\
&&-2k_2\varepsilon^2 \eta\big((1+\varepsilon\eta)u_0\big)_x-k_3\varepsilon\mu u_{0,xxx}+O(\varepsilon^2\mu,\mu^2)\nonumber \\
&=&- k_1\varepsilon\big((1+\varepsilon\eta)u_0\big)_x-2k_2\varepsilon^2(\eta u_0 )_x-2k_2\varepsilon^3u_0^2u_{0,x}-2k_2\varepsilon^3\eta (\eta u_0 )_x\nonumber \\
&&+(\frac{k_1}{6}-k_3)\varepsilon\mu u_{0,xxx}+O(\varepsilon^2\mu,\mu^2)\nonumber \\
&=&- k_1\varepsilon\big((1+(1+\frac{2k_2}{k_1})\varepsilon\eta)u_0\big)_x-2k_2\varepsilon^3(u^3_0)_x+(\frac{k_1}{6}-k_3)\varepsilon\mu u_{0,xxx}+O(\varepsilon^4,\varepsilon^2\mu,\mu^2).\nonumber \\
\end{eqnarray}
We impose the constraint
\begin{eqnarray}\label{relation1}
1+\frac{2k_2}{k_1}=k_1.
\end{eqnarray}
Thus, in view of the definition of $\varrho$ and (\ref{relationueta}), Eq. (\ref{varrhot}) becomes
\begin{eqnarray*}
{\varrho _t}=- k_1\varepsilon(\varrho u_0)_x+(2k_1k_2-2k_2)\varepsilon^3(u^3_0)_x\varrho+(\frac{k_1}{6}-k_3)\varepsilon\mu u_{0,xxx}+O(\varepsilon^4,\varepsilon^2\mu,\mu^2).
\end{eqnarray*}
To eliminate the dispersive term of $u_{0,xxx}$, with the choice of
\begin{eqnarray}\label{relation2}
\frac{k_1}{6}-k_3=0,
\end{eqnarray}
we have
\begin{eqnarray}\label{varrhot1}
{\varrho _t}+k_1\varepsilon(\varrho u_0)_x-(2k_1k_2-2k_2)\varepsilon^3(u^3_0)_x\varrho=O(\varepsilon^4,\varepsilon^2\mu,\mu^2).
\end{eqnarray}
Now we redefine the auxiliary function $\varrho$ in terms of the physical variable $\rho$ as follows
\begin{eqnarray}\label{rho}
\rho^{\frac{1}{b-1}}=\varrho, \quad b\neq1.
\end{eqnarray}
Substituting (\ref{rho}) into (\ref{varrhot1}), we obtain
\begin{eqnarray}\label{varrhotf}
{\rho _t}+k_1\varepsilon u_0 \rho_x+k_1\varepsilon(b-1)\rho u_{0,x}-(2k_1k_2-2k_2)(b-1)\varepsilon^3(u^3_0)_x\rho=O(\varepsilon^4,\varepsilon^2\mu,\mu^2).
\end{eqnarray}

On the other hand, it is inferred from (\ref{varrho}) retaining terms to order $O(\varepsilon^4,\varepsilon^2\mu)$ that
\begin{eqnarray*}\label{varrho2}
\varrho^2=1+2k_1\varepsilon \eta+2k_2\varepsilon^2u_0^2+(2k_2+k_1^2)\varepsilon^2 \eta^2+2k_1k_2\varepsilon^3(u_0^2+\eta^2)\eta+2k_3 \varepsilon\mu u_{0,xx}
+O(\varepsilon^4,\varepsilon^2\mu).
\end{eqnarray*}
Using the above equality, we express $\eta_x$ in terms of the auxiliary function $\varrho$ to find
\begin{eqnarray*}
\eta_x=\frac{1}{2k_1\varepsilon}\big(\varrho^2- 2k_2\varepsilon^2u_0^2- (2k_2+k_1^2)\varepsilon^2 \eta^2-
2k_1k_2\varepsilon^3(u_0^2+\eta^2)\eta-2k_3 \varepsilon\mu u_{0,xx}\big)_x+O(\varepsilon^3,\varepsilon\mu),
\end{eqnarray*}
which when substituted into (\ref{u0final}) yields
\begin{eqnarray}\label{etax1}
(u_{0}-\frac{\mu }{2}u_{0,xx})_t +(1-\frac{2k_2}{k_1}) \varepsilon {u_0}{u_{0,x}}&=&-\frac{1}{2k_1\varepsilon}(\varrho^2)_x+
\frac{2k_2+k_1^2}{2k_1}\varepsilon(\eta^2)_x+k_2\varepsilon^2\big((u_0^2+\eta^2)\eta\big)_x\nonumber \\
&&
+\frac{k_3}{k_1}\mu u_{0,xxx}+O(\varepsilon^3,\varepsilon\mu).
\end{eqnarray}
Using the relations (\ref{relationueta}) and the above equality (\ref{etax1}), and neglecting the terms of order $O(\varepsilon^3,\varepsilon\mu)$, we get
\begin{eqnarray}\label{etax2}
\big(u_{0}-(\frac{1}{2}-\frac{k_3}{k_1})\mu
u_{0,xx}\big)_t +(1-\frac{2k_2}{k_1}) \varepsilon {u_0}{u_{0,x}}&=&-\frac{1}{2k_1\varepsilon}(\varrho^2)_x+
\frac{2k_2+k_1^2}{2k_1}\varepsilon(\eta^2)_x+2k_2\varepsilon^2(u_0^3)_x.\nonumber \\
\end{eqnarray}
Now we reexpress the $\varrho^2$-term in (\ref{etax2}) in terms of $\rho^2$ by using
\begin{eqnarray}\label{rho1}
\varrho=1+\varepsilon\varrho_0,
\end{eqnarray}
where $\varrho_0$ is of order $O(1)$. Retaining terms to order $O(\varepsilon^4)$, we apply a Taylor expansion to write
\begin{eqnarray}\label{rho2}
\varrho^{2(b-1)}=(1+\varepsilon\varrho_0)^{2(b-1)}&=&1+2(b-1)\varepsilon\varrho_0+(b-1)(2b-3)\varepsilon^2\varrho^2_0\nonumber \\
&&+\frac{(b-1)(2b-3)(2b-4)}{3}\varepsilon^3\varrho^3_0+O(\varepsilon^4)=\rho^2,
\end{eqnarray}
and hence
\begin{eqnarray}\label{rho3}
2\varepsilon\varrho_0=\frac{1}{b-1}\rho^2-(2b-3)\varepsilon^2\varrho^2_0-\frac{(2b-3)(2b-4)}{3}\varepsilon^3\varrho^3_0-\frac{1}{b-1}+O(\varepsilon^4).
\end{eqnarray}
It then follows from (\ref{rho1}) that
\begin{eqnarray}\label{rho4}
\varrho^2=\frac{1}{b-1}\rho^2-(2b-4)\varepsilon^2\varrho^2_0-\frac{(2b-3)(2b-4)}{3}\varepsilon^3\varrho^3_0+\frac{b-2}{b-1}+O(\varepsilon^4).
\end{eqnarray}
Comparing (\ref{varrho}) and (\ref{rho1}), we find
\begin{eqnarray}\label{rho5}
\varrho_0=k_1\eta+k_2\varepsilon(u_0^2+\eta^2)+k_3\mu u_{0,xx},
\end{eqnarray}
and thus
\begin{eqnarray}\label{rho6}
\varrho_0^2=k^2_1\eta^2+2k_1k_2\varepsilon(u_0^2+\eta^2)\eta+O(\varepsilon^2,\mu),
\end{eqnarray}
and
\begin{eqnarray}\label{rho7}
\varrho^3_0=k^3_1\eta^3+O(\varepsilon).
\end{eqnarray}
Substituting (\ref{rho6})-(\ref{rho7}) into (\ref{rho4}), and using the relations (\ref{relationueta}),  we obtain
\begin{eqnarray}\label{rho8}
\varrho^2&=&\frac{1}{b-1}\rho^2-(2b-4)k_1^2\varepsilon^2\eta^2\nonumber \\
&&-(2b-4)k_1\cdot\frac{12k_2+(2b-3)k_1^2}{3}\varepsilon^3u_0^3+\frac{b-2}{b-1}+O(\varepsilon^4,\varepsilon^2\mu),
\end{eqnarray}
which allows us to replace the auxiliary function $\varrho^2$ appearing in (\ref{etax2}) with the physical variable $\rho^2$, in which case we deduce
\begin{eqnarray}\label{etax3}
&&\big(u_{0}-(\frac{1}{2}-\frac{k_3}{k_1})\mu
u_{0,xx}\big)_t +(1-\frac{2k_2}{k_1}) \varepsilon {u_0}{u_{0,x}}\nonumber \\
&=&-\frac{1}{2k_1\varepsilon}\cdot \frac{1}{b-1}(\rho^2)_x+
\frac{2k_2+k_1^2+(2b-4)k^2_1}{2k_1}\varepsilon(\eta^2)_x\nonumber \\
&&+\big(2k_2+(b-2)\cdot\frac{12k_2+(2b-3)k_1^2}{3}\big)\varepsilon^2(u_0^3)_x
+O(\varepsilon^3,\varepsilon\mu).
\end{eqnarray}
Note that we have only imposed two constraints (\ref{relation1})-(\ref{relation2})
on the coefficients $k_1$, $k_2$ and $k_3$. To eliminate the $(\eta^2)_x$-term,
we choose $k_1$ and $k_2$ such that
\begin{eqnarray}\label{relation3}
2k_2+k_1^2+(2b-4)k^2_1=0,
\end{eqnarray}
which upon substitution into (\ref{etax3}) gives
\begin{eqnarray}\label{etax4}
\varepsilon m_{0,t} +(1-\frac{2k_2}{k_1}) \varepsilon^2 {u_0}{u_{0,x}}
+\frac{1}{2k_1(b-1)}(\rho^2)_x+\beta\varepsilon^3(u_0^3)_x
=O(\varepsilon^4,\varepsilon^2\mu),
\end{eqnarray}
where $m_0:=u_{0}-(\frac{1}{2}-\frac{k_3}{k_1})\mu
u_{0,xx}=u_0-\frac{\mu}{3}u_{0,xx}$, and $\beta:=-2k_2-(b-2)\cdot\frac{12k_2+(2b-3)k_1^2}{3}$. Since $u_0=m_0+O(\mu)$, we have
\begin{eqnarray}\label{m}
{u_0}{u_{0,x}}=\frac{1}{1+b}(u_0m_{0,x}+b u_{0,x}m_0)+O(\mu).
\end{eqnarray}
Moreover, a simple computation shows that
\begin{eqnarray}\label{u3}
\beta\varepsilon^3(u_0^3)_x&=&\tilde{\beta}\varepsilon^3(u_0^3)_x+(\beta-\tilde{\beta})\varepsilon^3(u_0^3)_x\nonumber \\
&=&\tilde{\beta}\varepsilon^3(u_0^3\rho)_x+(\beta-\tilde{\beta})\varepsilon^3(u_0^3)_x+O(\varepsilon^4)\nonumber \\
&=&\frac{\tilde{\beta}}{k_1}\varepsilon^2(u_0^2\rho(k_1\varepsilon\eta))_x+(\beta-\tilde{\beta})\varepsilon^3(u_0^3)_x+O(\varepsilon^4)\nonumber \\
&=&\frac{\tilde{\beta}}{k_1(b-1)}\varepsilon^2(u_0^2\rho(\rho-1))_x+(\beta-\tilde{\beta})\varepsilon^3(u_0^3)_x+O(\varepsilon^4),
\end{eqnarray}
where $\tilde{\beta}$ is an arbitrary undetermined parameter. Combining (\ref{m}) and (\ref{u3}), Eq. (\ref{etax4}) can be rewritten as
\begin{eqnarray}\label{etax5}
&&\varepsilon m_{0,t} +(1-\frac{2k_2}{k_1})\cdot\frac{1}{1+b} \varepsilon^2 (u_0m_{0,x}+b u_{0,x}m_0)
+\frac{1}{2k_1(b-1)}(\rho^2)_x\nonumber \\
&&\quad\quad\quad\quad+\frac{\tilde{\beta}}{k_1(b-1)}\varepsilon^2(u_0^2\rho(\rho-1))_x+(\beta-\tilde{\beta})\varepsilon^3(u_0^3)_x
=O(\varepsilon^4,\varepsilon^2\mu,\mu^2).
\end{eqnarray}

Applying the following scaling
\begin{eqnarray}\label{scaling}
x\rightarrow \sqrt{\frac{\mu}{3}}x,\quad t\rightarrow \sqrt{\frac{\mu}{3}}t,\quad u_0\rightarrow \frac{1}{k_1 \varepsilon}u_0, \quad \rho\rightarrow \sqrt{b-1}\rho,
\end{eqnarray}
to (\ref{varrhotf}) and (\ref{etax5}), we establish a two-component $b$-family system of nonlinear equations which models two-dimensional shallow water waves of medium amplitude
\begin{equation}\label{system}
\left\{
 \begin{aligned}
 &\rho_t+u_0\rho_x+(b-1)\rho u_{0,x}-\frac{2k_1k_2-2k_2}{k_1^3}(b-1)(u_0^3)_x\rho=0,\\
&m_t+\frac{k_1-2k_2}{k_1^2(1+b)}(u_0m_x+b u_{0,x}m)+\frac{1}{2}(\rho^2)_x\\
&\quad\quad\quad\quad\quad\quad\quad\quad\quad\quad+\frac{\tilde{\beta}}{k_1^2}(u_0^2\rho^2)_x-\frac{\tilde{\beta}}{k_1^2\sqrt{b-1}}(u_0^2\rho)_x+\frac{\beta-\tilde{\beta}}{k_1^2}(u_0^3)_x=0,
\end{aligned}
\right.
\end{equation}
where the terms of order $O(\varepsilon^4,\varepsilon^2\mu,\mu^2)$ and higher have been neglected. Here, $m={u_{0}}- {u_{0,xx}},$ $b>1$ is a real number, and the constants $k_1$, $k_2$ and $k_3$
are completely determined by the constraints (\ref{relation1})-(\ref{relation2}) and (\ref{relation3}) in terms of $b$ as follows:
\begin{eqnarray}\label{coffe}
k_1=\frac{1}{2(b-1)},\quad k_2=-\frac{2b-3}{8(b-1)^2},\quad k_3=\frac{1}{12(b-1)},
\end{eqnarray}
and moreover,
\begin{eqnarray}\label{beta}
\beta=\frac{(2b-3)(5b-7)}{12(b-1)^2},
\end{eqnarray}
and the parameter $\tilde{\beta}$ is arbitrary and can be chosen appropriately as needed. For simplicity, we will drop the subscript and write $u$ instead of $u_0$ for the
system (\ref{system}) in the remainder of the paper.

\begin{remark2}\label{Rem2.1}
$(i)$ When $b=2$ and $\tilde{\beta}=-\frac{3}{4}$, the two-component $b$-family system (\ref{system}) reduces exactly to system (\ref{2CH}), i.e., the two-component Camassa-Holm system under the CH scaling.\\
$(ii)$ In the medium-amplitude regime ($\mu\ll1$ and $\varepsilon=O(\sqrt{\mu})$), the unscaled physical variable
$\rho$ satisfies $\rho=1+O(\varepsilon)$. After the scaling (\ref{scaling}), this relation becomes
$\rho=(b-1)^{-1/2}+O(\varepsilon)$. These facts allow us to simplify (\ref{system}) asymptotically.
Indeed, replacing the unscaled variable $\rho$ by $1$ in the last term of (\ref{varrhotf}) introduces an error
of $O(\varepsilon^4)$, since that term already contains a factor $\varepsilon^3$. Similarly, in (\ref{u3}),
the identity $\rho(\rho-1)=\rho-1+O(\varepsilon^2)$ allows us to replace $u_0^2\rho(\rho-1)$ by
$u_0^2(\rho-1)$ with an error of $O(\varepsilon^4)$. Both errors are of the same order as the terms already
neglected in the derivation. After these simplifications and the scaling (\ref{scaling}), we obtain the
following asymptotically equivalent system:
\begin{equation}\label{system1}
\left\{
 \begin{aligned}
 &\rho_t+u\rho_x+(b-1)\rho u_{x}-\frac{2k_1k_2-2k_2}{k_1^3}\sqrt{b-1}(u^3)_x=0,\\
&m_t+\frac{k_1-2k_2}{k_1^2(1+b)}(um_x+b u_{x}m)+\frac{1}{2}(\rho^2)_x\\
&\quad\quad\quad\quad\quad\quad\quad\quad\quad\quad+\frac{\tilde{\beta}}{k_1^2\sqrt{b-1}}(u^2\rho)_x-\frac{\tilde{\beta}}{k_1^2(b-1)}(u^2)_x+\frac{\beta-\tilde{\beta}}{k_1^2}(u^3)_x=0,
\end{aligned}
\right.
\end{equation}
where $m=u-u_{xx}$, the coefficients $k_1$, $k_2$, $\beta$ are given by (\ref{coffe})-(\ref{beta}), and $\tilde{\beta}$ is still a free parameter. This model (\ref{system1}) preserves the essential nonlinearity of
medium-amplitude waves while being mathematically simpler. In particular, taking $b=2$ and $\tilde{\beta}=-\frac{3}{4}$ in (\ref{system1}), we obtain the following model
\begin{equation}\label{systems}
\left\{
 \begin{aligned}
 &\rho_t+(\rho u)_x-(u^3)_x=0,\\
&m_t+(um_x+2 u_{x}m)+\frac{1}{2}(\rho^2)_x-3(u^2(\rho-1))_x+4(u^3)_x=0.
\end{aligned}
\right.
\end{equation}
where $m=u-u_{xx}$. We call this system (\ref{systems}) the simplified two-component Camassa-Holm system in the medium-amplitude regime, denoted as the s2CH system for brevity.
\end{remark2}

\section{Local well-posedness}\label{sec3}
\newtheorem{theorem3}{Theorem}[section]
\newtheorem{lemma3}{Lemma}[section]
\newtheorem {remark3}{Remark}[section]
\newtheorem {definition3}{Definition}[section]
\newtheorem{corollary3}{Corollary}[section]
\newtheorem{example3}{Example}[section]
\par
In this section, we establish the local well-posedness of the Cauchy problem for the system (\ref{1.1}) in the
nonhomogeneous Besov spaces by using the Littlewood-Paley theory and the transport equation theory.  Since all function spaces in this and the next
sections are over $\R$, for simplicity, we drop $\R$ in the notation of function spaces if there is no ambiguity.

Denote by $G(x)=\frac12e^{-|x|}$, $x\in\R$, the kernel of the operator $(1-\partial^2_x)^{-1}$, so that
$G\ast m=u$ for $m=u-u_{xx}$. We can rewrite the system (\ref{1.1}) in the following transport-nonlocal form:
\begin{equation}\label{nonlocal}
\left\{
 \begin{aligned}
 &\rho_t+u\rho_x=\big(3\gamma_1u^2-(b-1)\big)\rho u_x,\\
&u_t+\alpha u u_x=-\partial_xG\ast F(u,\rho),
\end{aligned}
\right.
\end{equation}
where
\begin{equation}\label{F}
F(u,\rho)=\frac{\alpha b}{2}u^2+\frac{\alpha(3-b)}{2}u_x^2+\frac12\rho^2+\gamma_2u^2\rho^2+\gamma_3u^2\rho+\gamma_4u^3.
\end{equation}

Our first result of this section reads as follows.
\begin{theorem3}\label{LWP}
Let $1\leq p,r\leq\infty$, $s>\max\{1+\frac{1}{p},\frac{3}{2}\}$ and $(\rho_0,u_0)\in B^{s-1}_{p,r}\times B^{s}_{p,r}$.
Then there exists some $T>0$, depending only on $\|(\rho_0,u_0)\|_{B^{s-1}_{p,r}\times B^{s}_{p,r}}$, such that the
system (\ref{1.1}) has a unique solution $(\rho,u)$ in
\begin{equation}\label{Es}
E^s_{p,r}(T):=
\left\{
 \begin{aligned}
&C([0,T);B^{s-1}_{p,r})\times C([0,T);B^{s}_{p,r})\cap C^1([0,T);B^{s-2}_{p,r})\times C^1([0,T);B^{s-1}_{p,r}),\\
 &if\ r<\infty,\\
&C_w([0,T);B^{s-1}_{p,\infty})\times C_w([0,T);B^{s}_{p,\infty})\cap C^{0,1}([0,T);B^{s-2}_{p,\infty})\times
C^{0,1}([0,T);B^{s-1}_{p,\infty}),\\
&if\ r=\infty.
\end{aligned}
\right.
\end{equation}
Moreover, the solution map $(\rho_0,u_0)\mapsto(\rho,u)$ is H\"older continuous from
$B^{s-1}_{p,r}\times B^{s}_{p,r}$ to $E^{s'}_{p,r}(T)$ with $s-1\leq s'<s$, that is,
\begin{equation}\label{holder}
\|\rho_1-\rho_2\|_{L^\infty(0,T;B^{s'-1}_{p,r})}+\|u_1-u_2\|_{L^\infty(0,T;B^{s'}_{p,r})}
\leq C\big(\|\rho_1(0)-\rho_2(0)\|_{B^{s-1}_{p,r}}^{\theta}+\|u_1(0)-u_2(0)\|_{B^{s}_{p,r}}^{\theta}\big),
\end{equation}
where $\theta=s-s'\in(0,1]$.
\end{theorem3}
\begin{proof}
In order to prove Theorem \ref{LWP}, we proceed in four steps.\\
\textbf{Step 1: Construction of approximate solutions.}
Starting from $(\rho_1,u_1):=(S_1\rho_0,S_1u_0)$, we define by induction sequences $(\rho_n,u_n)_{n\in\mathbb{N}}$
by solving the following linear transport equations:
\begin{equation}\label{appro}
\left\{
 \begin{aligned}
 &\partial_t\rho_{n+1}+u_n\partial_x\rho_{n+1}=\big(3\gamma_1u_n^2-(b-1)\big)\rho_n\partial_xu_n,\\
&\partial_tu_{n+1}+\alpha u_n\partial_xu_{n+1}=-\partial_xG\ast F(u_n,\rho_n),\\
&(\rho_{n+1},u_{n+1})|_{t=0}=(S_{n+1}\rho_0,S_{n+1}u_0),
\end{aligned}
\right.
\end{equation}
where $F(u,\rho)$ is given by (\ref{F}) and $S_{n+1}$ is the low-frequency cut-off operator defined in the Appendix.
Assume that $(\rho_n,u_n)\in L^\infty(0,T;B^{s-1}_{p,r})\times L^\infty(0,T;B^{s}_{p,r})$. Since
$s>\max\{1+\frac{1}{p},\frac{3}{2}\}$, we have $s-1>\max\{\frac1p,\frac12\}$ and hence $B^{s-1}_{p,r}$ is an
algebra, which implies
$\big(3\gamma_1u_n^2-(b-1)\big)\rho_n\partial_xu_n\in L^\infty(0,T;B^{s-1}_{p,r})$. Moreover, the operator
$\partial_xG\ast=\partial_x(1-\partial^2_x)^{-1}$ is a multiplier of degree $-1$, which together with Corollary \ref{algebra} yields
\begin{eqnarray}\label{Fest}
\|\partial_xG\ast F(u_n,\rho_n)\|_{B^s_{p,r}}&\leq& C\|F(u_n,\rho_n)\|_{B^{s-1}_{p,r}}\nonumber\\
&\leq&C\Big(\|u_n\|^2_{B^s_{p,r}}+\|\partial_xu_n\|^2_{B^{s-1}_{p,r}}+\|\rho_n\|^2_{B^{s-1}_{p,r}}
+\|u_n\|^2_{B^{s}_{p,r}}\|\rho_n\|^2_{B^{s-1}_{p,r}}\nonumber\\
&&\quad+\|u_n\|^2_{B^{s}_{p,r}}\|\rho_n\|_{B^{s-1}_{p,r}}+\|u_n\|^3_{B^s_{p,r}}\Big)\nonumber\\
&\leq&C\Big(N_n^2(t)+N_n^4(t)\Big),
\end{eqnarray}
where we denote
\begin{equation*}
N_n(t):=\sup_{\tau\in[0,t]}\Big(\|\rho_n(\tau)\|_{B^{s-1}_{p,r}}+\|u_n(\tau)\|_{B^{s}_{p,r}}\Big).
\end{equation*}
Hence $-\partial_xG\ast F(u_n,\rho_n)\in L^\infty(0,T;B^{s}_{p,r})$. The transport equation theory
(see Chapter 3 in \cite{BCD}) then ensures that
$(\rho_{n+1},u_{n+1})\in L^\infty(0,T;B^{s-1}_{p,r})\times L^\infty(0,T;B^{s}_{p,r})$.\\
\textbf{Step 2: Uniform bounds for the approximate solutions.} We now find some positive $T$ such that the approximate solutions are uniformly bounded on $[0,T]$.
Assume first $s\neq2+\frac1p$. Applying Lemma \ref{transport} with $p_1=p$ to the first equation of (\ref{appro})
(note that $V'(t)\leq C\|u_n(t)\|_{B^s_{p,r}}$ in both cases $s-1\leq1+\frac1p$ and $s-1>1+\frac1p$, since
$s>\max\{1+\frac1p,\frac32\}$), we infer that
\begin{eqnarray}\label{rhoest}
\|\rho_{n+1}\|_{L^\infty_t(B^{s-1}_{p,r})}&\leq& e^{CU_n(t)}\Big(\|S_{n+1}\rho_0\|_{B^{s-1}_{p,r}}
+\int_0^te^{-CU_n(\tau)}\big\|\big(3\gamma_1u_n^2-(b-1)\big)\rho_n\partial_xu_n\big\|_{B^{s-1}_{p,r}}d\tau\Big)\nonumber\\
&\leq&e^{CU_n(t)}\Big(\|\rho_0\|_{B^{s-1}_{p,r}}
+C\int_0^te^{-CU_n(\tau)}\big(1+\|u_n\|^2_{B^s_{p,r}}\big)\|\rho_n\|_{B^{s-1}_{p,r}}\|u_n\|_{B^s_{p,r}}d\tau\Big),\nonumber\\
\end{eqnarray}
where $U_n(t):=\int_0^t\|u_n(\tau)\|_{B^s_{p,r}}d\tau$. Similarly, applying Lemma \ref{transport} to the second
equation of (\ref{appro}), and using (\ref{Fest}), we obtain
\begin{eqnarray}\label{uest}
\|u_{n+1}\|_{L^\infty_t(B^{s}_{p,r})}\leq e^{CU_n(t)}\Big(\|u_0\|_{B^{s}_{p,r}}
+C\int_0^te^{-CU_n(\tau)}\big(N_n^2(\tau)+N_n^4(\tau)\big)d\tau\Big).
\end{eqnarray}
Adding (\ref{rhoest}) and (\ref{uest}), we deduce that
\begin{eqnarray}\label{Nest}
N_{n+1}(t)\leq e^{CU_n(t)}\Big(N_0+C\int_0^t\big(N_n^2(\tau)+N_n^4(\tau)\big)d\tau\Big),
\end{eqnarray}
where $N_0:=\|\rho_0\|_{B^{s-1}_{p,r}}+\|u_0\|_{B^{s}_{p,r}}$. Now we fix $T>0$ such that
\begin{equation}\label{Tcond}
2CN_0T\leq\ln\frac32\quad\quad\mbox{and}\quad\quad 4CTN_0\big(1+4N_0^2\big)\leq\frac13,
\end{equation}
and suppose that
\begin{equation}\label{induction}
\forall t\in[0,T],\quad N_n(t)\leq2N_0.
\end{equation}
Then $U_n(t)\leq2N_0T$, which together with (\ref{Nest}) implies that
\begin{eqnarray*}
N_{n+1}(t)\leq e^{2CN_0T}\Big(N_0+4CTN_0^2\big(1+4N_0^2\big)\Big)
\leq\frac32\Big(N_0+\frac{N_0}{3}\Big)=2N_0.
\end{eqnarray*}
Therefore, by induction, we obtain that $(\rho_n,u_n)$ is bounded in
$L^\infty(0,T;B^{s-1}_{p,r})\times L^\infty(0,T;B^{s}_{p,r})$. In the endpoint case $s=2+\frac1p$, the same
conclusion follows by applying Lemma \ref{transportcri} instead of Lemma \ref{transport} and repeating the above
argument, which shows that $(\rho_n,u_n)$ is bounded in
$L^\infty(0,T;B^{1+\frac1p}_{p,r})\times L^\infty(0,T;B^{2+\frac1p}_{p,r})$.\\
\textbf{Step 3: Convergence of the approximate sequence.}
 We are going to show that $\{(\rho_n,u_n)\}$ is a Cauchy sequence in
$L^\infty(0,T;B^{s-2}_{p,r})\times L^\infty(0,T;B^{s-1}_{p,r})$. For this purpose, we deduce from (\ref{appro}) that
\begin{equation}\label{diff}
\left\{
 \begin{aligned}
 &\partial_t(\rho_{n+m+1}-\rho_{n+1})+u_{n+m}\partial_x(\rho_{n+m+1}-\rho_{n+1})
=(u_n-u_{n+m})\partial_x\rho_{n+1}+R^1_{n,m},\\
&\partial_t(u_{n+m+1}-u_{n+1})+\alpha u_{n+m}\partial_x(u_{n+m+1}-u_{n+1})
=\alpha(u_n-u_{n+m})\partial_xu_{n+1}+R^2_{n,m},
\end{aligned}
\right.
\end{equation}
where
\begin{eqnarray*}
R^1_{n,m}&=&\big(3\gamma_1u_{n+m}^2-(b-1)\big)\rho_{n+m}\partial_xu_{n+m}
-\big(3\gamma_1u_{n}^2-(b-1)\big)\rho_{n}\partial_xu_{n},\\
R^2_{n,m}&=&-\partial_xG\ast\big(F(u_{n+m},\rho_{n+m})-F(u_{n},\rho_{n})\big).
\end{eqnarray*}
By virtue of Lemma \ref{transport} (and Lemma \ref{transportcri} in the endpoint case $s=2+\frac1p$)
and using the fact that $\{(\rho_n,u_n)\}$ is bounded in
$L^\infty(0,T;B^{s-1}_{p,r})\times L^\infty(0,T;B^{s}_{p,r})$, we infer that
\begin{eqnarray}\label{diff1}
&&\|\rho_{n+m+1}-\rho_{n+1}\|_{L^\infty(0,T;B^{s-2}_{p,r})}\nonumber\\
&\leq&C\Big(\|S_{n+m+1}\rho_0-S_{n+1}\rho_0\|_{B^{s-2}_{p,r}}
+\int_0^T\|(u_n-u_{n+m})\partial_x\rho_{n+1}\|_{B^{s-2}_{p,r}}+\|R^1_{n,m}\|_{B^{s-2}_{p,r}}d\tau\Big),\nonumber\\
\end{eqnarray}
and
\begin{eqnarray}\label{diff2}
&&\|u_{n+m+1}-u_{n+1}\|_{L^\infty(0,T;B^{s-1}_{p,r})}\nonumber\\
&\leq&C\Big(\|S_{n+m+1}u_0-S_{n+1}u_0\|_{B^{s-1}_{p,r}}
+\int_0^T\|(u_n-u_{n+m})\partial_xu_{n+1}\|_{B^{s-1}_{p,r}}+\|R^2_{n,m}\|_{B^{s-1}_{p,r}}d\tau\Big).\nonumber\\
\end{eqnarray}
Taking advantage of Lemma \ref{Morse} with $\sigma=s-1$ (note that $\sigma=s-1>\max\{\frac1p,\frac12\}$), we have
\begin{eqnarray}\label{diff3}
&&\|(u_n-u_{n+m})\partial_x\rho_{n+1}\|_{B^{s-2}_{p,r}}
\leq C\|u_n-u_{n+m}\|_{B^{s-1}_{p,r}}\|\partial_x\rho_{n+1}\|_{B^{s-2}_{p,r}}\nonumber\\
&\leq&C\|u_n-u_{n+m}\|_{B^{s-1}_{p,r}}\|\rho_{n+1}\|_{B^{s-1}_{p,r}}
\leq C\|u_n-u_{n+m}\|_{B^{s-1}_{p,r}}.
\end{eqnarray}
For $R^1_{n,m}$, we decompose it as
\begin{eqnarray*}
R^1_{n,m}&=&3\gamma_1(u_{n+m}^2-u_{n}^2)\rho_{n+m}\partial_xu_{n+m}
+3\gamma_1u_{n}^2(\rho_{n+m}-\rho_{n})\partial_xu_{n+m}\\
&&+3\gamma_1u_{n}^2\rho_{n}\partial_x(u_{n+m}-u_{n})
-(b-1)(\rho_{n+m}-\rho_{n})\partial_xu_{n+m}-(b-1)\rho_{n}\partial_x(u_{n+m}-u_{n}).
\end{eqnarray*}
Since $B^{s-1}_{p,r}$ is an algebra, by Lemma \ref{Morse} and the uniform boundedness of $\{(\rho_n,u_n)\}$, we get
\begin{eqnarray}\label{R1est}
\|R^1_{n,m}\|_{B^{s-2}_{p,r}}&\leq&C\Big(\|u_{n+m}^2-u_{n}^2\|_{B^{s-1}_{p,r}}
\|\rho_{n+m}\partial_xu_{n+m}\|_{B^{s-1}_{p,r}}
+\|\rho_{n+m}-\rho_{n}\|_{B^{s-2}_{p,r}}\|u_{n}^2\partial_xu_{n+m}\|_{B^{s-1}_{p,r}}\nonumber\\
&&\quad+\|\partial_x(u_{n+m}-u_{n})\|_{B^{s-2}_{p,r}}\|u_{n}^2\rho_{n}\|_{B^{s-1}_{p,r}}
+\|\rho_{n+m}-\rho_{n}\|_{B^{s-2}_{p,r}}\|\partial_xu_{n+m}\|_{B^{s-1}_{p,r}}\nonumber\\
&&\quad+\|\rho_{n}\|_{B^{s-1}_{p,r}}\|\partial_x(u_{n+m}-u_{n})\|_{B^{s-2}_{p,r}}\Big)\nonumber\\
&\leq&C\Big(\|u_{n+m}-u_{n}\|_{B^{s-1}_{p,r}}+\|\rho_{n+m}-\rho_{n}\|_{B^{s-2}_{p,r}}\Big).
\end{eqnarray}
Plugging (\ref{diff3}) and (\ref{R1est}) into (\ref{diff1}) yields that
\begin{eqnarray}\label{rhoCauchy}
&&\|\rho_{n+m+1}-\rho_{n+1}\|_{L^\infty(0,T;B^{s-2}_{p,r})}\nonumber\\
&\leq& C\Big(\|S_{n+m+1}\rho_0-S_{n+1}\rho_0\|_{B^{s-2}_{p,r}}
+\int_0^T\|u_{n+m}-u_{n}\|_{B^{s-1}_{p,r}}+\|\rho_{n+m}-\rho_{n}\|_{B^{s-2}_{p,r}}d\tau\Big).\nonumber\\
\end{eqnarray}
Similarly, since $B^{s-1}_{p,r}$ is an algebra, we deduce that
\begin{eqnarray}\label{diff4}
\|(u_n-u_{n+m})\partial_xu_{n+1}\|_{B^{s-1}_{p,r}}
\leq C\|u_n-u_{n+m}\|_{B^{s-1}_{p,r}}\|\partial_xu_{n+1}\|_{B^{s-1}_{p,r}}
\leq C\|u_n-u_{n+m}\|_{B^{s-1}_{p,r}}.\nonumber\\
\end{eqnarray}
It remains to estimate $R^2_{n,m}$. Note that $\|\partial_xG\ast f\|_{B^{s-1}_{p,r}}\leq C\|f\|_{B^{s-2}_{p,r}}$, so we only need to bound
$\|F(u_{n+m},\rho_{n+m})-F(u_{n},\rho_{n})\|_{B^{s-2}_{p,r}}$. Each term is treated by adding and subtracting
suitable terms and then applying Lemma \ref{Morse} and the uniform bound. For instance,
\begin{eqnarray}\label{d1}
\|u_{n+m}^2-u_{n}^2\|_{B^{s-2}_{p,r}}
&\leq&\|(u_{n+m}-u_{n})(u_{n+m}+u_{n})\|_{B^{s-1}_{p,r}}
\leq C\|u_{n+m}-u_{n}\|_{B^{s-1}_{p,r}},
\end{eqnarray}
\begin{eqnarray}\label{d2}
\|(\partial_xu_{n+m})^2-(\partial_xu_{n})^2\|_{B^{s-2}_{p,r}}
&\leq& C\|\partial_x(u_{n+m}-u_{n})\|_{B^{s-2}_{p,r}}\|\partial_x(u_{n+m}+u_{n})\|_{B^{s-1}_{p,r}}\nonumber\\
&\leq& C\|u_{n+m}-u_{n}\|_{B^{s-1}_{p,r}},
\end{eqnarray}
and for the quartic term,
\begin{eqnarray}\label{d4}
\|u_{n+m}^2\rho_{n+m}^2-u_{n}^2\rho_{n}^2\|_{B^{s-2}_{p,r}}
&\leq&C\Big(\|u_{n+m}^2-u_{n}^2\|_{B^{s-2}_{p,r}}\|\rho_{n+m}^2\|_{B^{s-1}_{p,r}}
+\|u_{n}^2\|_{B^{s-1}_{p,r}}\|\rho_{n+m}^2-\rho_{n}^2\|_{B^{s-2}_{p,r}}\Big)\nonumber\\
&\leq&C\Big(\|u_{n+m}-u_{n}\|_{B^{s-1}_{p,r}}+\|\rho_{n+m}-\rho_{n}\|_{B^{s-2}_{p,r}}\Big),
\end{eqnarray}
where in the last step we used the decomposition
$u_{n+m}^2\rho_{n+m}^2-u_{n}^2\rho_{n}^2=(u_{n+m}^2-u_{n}^2)\rho_{n+m}^2+u_{n}^2(\rho_{n+m}^2-\rho_{n}^2)$.
The remaining terms $\|\rho_{n+m}^2-\rho_{n}^2\|_{B^{s-2}_{p,r}}$,
$\|u_{n+m}^2\rho_{n+m}-u_{n}^2\rho_{n}\|_{B^{s-2}_{p,r}}$ and
$\|u_{n+m}^3-u_{n}^3\|_{B^{s-2}_{p,r}}$ can be estimated in exactly the same way. Collecting these estimates, we obtain
\begin{eqnarray}\label{R2est}
\|R^2_{n,m}\|_{B^{s-1}_{p,r}}
\leq C\Big(\|u_{n+m}-u_{n}\|_{B^{s-1}_{p,r}}+\|\rho_{n+m}-\rho_{n}\|_{B^{s-2}_{p,r}}\Big),
\end{eqnarray}
which along with (\ref{diff2}) and (\ref{diff4}) leads to
\begin{eqnarray}\label{uCauchy}
&&\|u_{n+m+1}-u_{n+1}\|_{L^\infty(0,T;B^{s-1}_{p,r})}\nonumber\\
&\leq& C\Big(\|S_{n+m+1}u_0-S_{n+1}u_0\|_{B^{s-1}_{p,r}}
+\int_0^T\|u_{n+m}-u_{n}\|_{B^{s-1}_{p,r}}+\|\rho_{n+m}-\rho_{n}\|_{B^{s-2}_{p,r}}d\tau\Big).\nonumber\\
\end{eqnarray}
The inequality (\ref{uCauchy}) together with (\ref{rhoCauchy}) yields that
\begin{eqnarray}\label{Cauchy}
&&\|\rho_{n+m+1}-\rho_{n+1}\|_{L^\infty(0,T;B^{s-2}_{p,r})}
+\|u_{n+m+1}-u_{n+1}\|_{L^\infty(0,T;B^{s-1}_{p,r})}\nonumber\\
&\leq& C\Big(\|S_{n+m+1}\rho_0-S_{n+1}\rho_0\|_{B^{s-1}_{p,r}}+\|S_{n+m+1}u_0-S_{n+1}u_0\|_{B^{s}_{p,r}}\nonumber\\
&&\quad+\int_0^T\|u_{n+m}-u_{n}\|_{B^{s-1}_{p,r}}+\|\rho_{n+m}-\rho_{n}\|_{B^{s-2}_{p,r}}d\tau\Big).
\end{eqnarray}
By defining $A_{n,m}(t):=\|\rho_{n+m}-\rho_{n}\|_{B^{s-2}_{p,r}}+\|u_{n+m}-u_{n}\|_{B^{s-1}_{p,r}}$,
we deduce that
\begin{eqnarray}\label{Arec}
A_{n+1,m}(t)\leq C\Big(A_{n+1,m}(0)+\int_0^tA_{n,m}(\tau)d\tau\Big),\qquad n\geq1.
\end{eqnarray}
Arguing by induction with respect to $n$ in (\ref{Arec}), we get, for every $n\geq1$,
\begin{eqnarray}\label{Arec2}
\sup_{t\in[0,T]}A_{n+1,m}(t)\leq\sum_{k=0}^{n-1} C\frac{(CT)^k}{k!}A_{n+1-k,m}(0)
+C\frac{(CT)^n}{n!},
\end{eqnarray}
where the uniform bound for $A_{1,m}$ has been absorbed into the last term. This leads to $A_{n,m}(t)\rightarrow0$ as
$n\rightarrow\infty$, $\forall m\in\mathbb{N}$. So $\{(\rho_n,u_n)\}$ is a Cauchy sequence in
$L^\infty(0,T;B^{s-2}_{p,r})\times L^\infty(0,T;B^{s-1}_{p,r})$ and converges to some limit function
$(\rho,u)\in L^\infty(0,T;B^{s-2}_{p,r})\times L^\infty(0,T;B^{s-1}_{p,r})$. Since $(\rho_n,u_n)$ is bounded in
$L^\infty(0,T;B^{s-1}_{p,r})\times L^\infty(0,T;B^{s}_{p,r})$, the Fatou property for Besov spaces
(see Remark \ref{Fatou}) yields that $(\rho,u)$ also belongs to
$L^\infty(0,T;B^{s-1}_{p,r})\times L^\infty(0,T;B^{s}_{p,r})$. Passing to the limit in (\ref{appro}), we deduce
that $(\rho,u)$ is indeed a solution of (\ref{1.1}). By the transport equations
theory \cite{BCD}, we infer that $(\rho,u)\in E^s_{p,r}(T)$.\\
\textbf{Step 4: Uniqueness and stability.} Suppose that $(\rho_1,u_1)$ and
$(\rho_2,u_2)$ are two solutions of (\ref{1.1}) with the initial data $(\rho_1(0),u_1(0))$ and $(\rho_2(0),u_2(0))$,
respectively. Then their difference satisfies
\begin{equation}\label{uniq}
\left\{
 \begin{aligned}
 &\partial_t(\rho_1-\rho_2)+u_1\partial_x(\rho_1-\rho_2)
=(u_2-u_1)\partial_x\rho_2+\big(3\gamma_1u_1^2-(b-1)\big)\rho_1\partial_xu_1\\
&\quad\quad\quad\quad\quad\quad\quad\quad\quad\quad\quad\quad\quad-\big(3\gamma_1u_2^2-(b-1)\big)\rho_2\partial_xu_2,\\
&\partial_t(u_1-u_2)+\alpha u_1\partial_x(u_1-u_2)
=\alpha(u_2-u_1)\partial_xu_2-\partial_xG\ast\big(F(u_1,\rho_1)-F(u_2,\rho_2)\big).
\end{aligned}
\right.
\end{equation}
Applying the transport estimates used in Step 3 to (\ref{uniq}) and estimating all the source terms by the
same calculations, we arrive at
\begin{eqnarray}\label{uniq6}
&&\|\rho_1(t)-\rho_2(t)\|_{B^{s-2}_{p,r}}+\|u_1(t)-u_2(t)\|_{B^{s-1}_{p,r}}
\leq C\Big(\|\rho_1(0)-\rho_2(0)\|_{B^{s-1}_{p,r}}+\|u_1(0)-u_2(0)\|_{B^{s}_{p,r}}
\nonumber\\
&&\quad +\int_0^t\|\rho_1(\tau)-\rho_2(\tau)\|_{B^{s-2}_{p,r}}+\|u_1(\tau)-u_2(\tau)\|_{B^{s-1}_{p,r}}d\tau\Big).\nonumber
\end{eqnarray}
Using Gronwall's inequality, we get
\begin{eqnarray}\label{uniq7}
&&\sup_{t\in[0,T)}\|\rho_1(t)-\rho_2(t)\|_{B^{s-2}_{p,r}}+\sup_{t\in[0,T)}\|u_1(t)-u_2(t)\|_{B^{s-1}_{p,r}}\nonumber\\
&\leq& C\Big(\|\rho_1(0)-\rho_2(0)\|_{B^{s-1}_{p,r}}+\|u_1(0)-u_2(0)\|_{B^{s}_{p,r}}\Big).
\end{eqnarray}
Taking advantage of the interpolation argument ensures that
\begin{eqnarray*}\label{uniq8}
&&\sup_{t\in[0,T]}\|\rho_1(t)-\rho_2(t)\|_{B^{s'-1}_{p,r}}+\sup_{t\in[0,T]}\|u_1(t)-u_2(t)\|_{B^{s'}_{p,r}}\nonumber\\
&\leq &C\Big(\|\rho_1(0)-\rho_2(0)\|_{B^{s-1}_{p,r}}^{\theta}+\|u_1(0)-u_2(0)\|_{B^{s}_{p,r}}^{\theta}\Big),
\end{eqnarray*}
where $\theta=s-s'\in(0,1]$. This is exactly (\ref{holder}).
In particular, taking $(\rho_1(0),u_1(0))=(\rho_2(0),u_2(0))$ in (\ref{uniq7}) yields the uniqueness.
This completes the proof of Theorem \ref{LWP}.
\end{proof}

Next, we establish the local
well-posedness of (\ref{1.1}) in the critical Besov space $B^{\frac12}_{2,1}\times B^{\frac32}_{2,1}$, which is
the second main result of this section.
\begin{theorem3}\label{LWPcri}
Let $(\rho_0,u_0)\in B^{\frac12}_{2,1}\times B^{\frac32}_{2,1}$. Then there exists some $T>0$, depending only on
$\|(\rho_0,u_0)\|_{B^{\frac12}_{2,1}\times B^{\frac32}_{2,1}}$, such that the system (\ref{1.1}) has a unique
solution $(\rho,u)$ in
\begin{equation*}
C([0,T);B^{\frac12}_{2,1})\times C([0,T);B^{\frac32}_{2,1})\cap C^1([0,T);B^{-\frac12}_{2,1})\times
C^1([0,T);B^{\frac12}_{2,1}).
\end{equation*}
Moreover, for any $s'\in[\frac12,\frac32)$, there holds
\begin{eqnarray}\label{holdercri}
&&\sup_{t\in[0,T)}\|\rho_1(t)-\rho_2(t)\|_{B^{s'-1}_{2,1}}+\sup_{t\in[0,T)}
\|u_1(t)-u_2(t)\|_{B^{s'}_{2,1}}\nonumber\\
&\leq & C\Big(\|\rho_1(0)-\rho_2(0)\|_{B^{\frac12}_{2,1}}+\|u_1(0)-u_2(0)\|_{B^{\frac32}_{2,1}}\Big)^{\theta\exp\{-CT\}},
\end{eqnarray}
where $\theta=\frac32-s'\in(0,1]$.
\end{theorem3}
\begin{proof}
We still use the approximate sequence $(\rho_n,u_n)$ constructed in (\ref{appro}) and proceed in four steps.\\
\textbf{Step 1: Construction of approximate solutions.}
Since $B^{\frac12}_{2,1}$ is an algebra (see Corollary \ref{algebra}) and
$B^{\frac32}_{2,1}\hookrightarrow B^{\frac12}_{2,1}\hookrightarrow L^\infty$, one checks as in Step 1 of the proof
of Theorem \ref{LWP} that the right-hand sides of (\ref{appro}) belong to $L^\infty(0,T;B^{\frac12}_{2,1})$ and
$L^\infty(0,T;B^{\frac32}_{2,1})$, respectively. The theory of transport equations then ensures that the iterates
are well-defined.\\
\textbf{Step 2: Uniform bounds for the approximate solutions.}
For both equations of (\ref{appro}) we have $V'(t)\leq C\|u_n(t)\|_{B^{\frac32}_{2,1}}$.
Hence the argument of Step 2 in the proof of Theorem \ref{LWP} applies similarly with
$N_n(t):=\sup_{\tau\in[0,t]}\big(\|\rho_n(\tau)\|_{B^{\frac12}_{2,1}}+\|u_n(\tau)\|_{B^{\frac32}_{2,1}}\big)$,
and yields some $T>0$ as in (\ref{Tcond}) such that
\begin{equation}\label{cribound}
N_n(t)\leq2N_0,\qquad \forall t\in[0,T],\ \forall n\in\mathbb{N},
\end{equation}
where $N_0:=\|\rho_0\|_{B^{\frac12}_{2,1}}+\|u_0\|_{B^{\frac32}_{2,1}}$.\\
\textbf{Step 3: Convergence of the approximate sequence.}
We are going to show that $\{(\rho_n,u_n)\}$ is a Cauchy sequence in
$L^\infty(0,T;B^{-\frac12}_{2,\infty})\times L^\infty(0,T;B^{\frac12}_{2,\infty})$. Applying Lemma \ref{transport}
to (\ref{diff}), we claim that
\begin{eqnarray}\label{cridiff}
A_{n+1,m}(t)\leq C\Big(a_{n,m}+\int_0^t\big(\|\rho_{n+m}-\rho_n\|_{B^{-\frac12}_{2,\infty}}
+\|u_{n+m}-u_n\|_{B^{\frac12}_{2,1}}\big)d\tau\Big),
\end{eqnarray}
where
\begin{eqnarray*}
A_{n,m}(t)&:=&\|\rho_{n+m}(t)-\rho_n(t)\|_{B^{-\frac12}_{2,\infty}}
+\|u_{n+m}(t)-u_n(t)\|_{B^{\frac12}_{2,\infty}},\\
a_{n,m}&:=&\|S_{n+m+1}\rho_0-S_{n+1}\rho_0\|_{B^{-\frac12}_{2,\infty}}
+\|S_{n+m+1}u_0-S_{n+1}u_0\|_{B^{\frac12}_{2,\infty}}.
\end{eqnarray*}
Indeed, all the source terms in (\ref{diff}) can be estimated by Lemma \ref{criproduct} together with the algebra
property of $B^{\frac12}_{2,1}$ and (\ref{cribound}). First,
\begin{eqnarray}\label{crid1}
\|(u_n-u_{n+m})\partial_x\rho_{n+1}\|_{B^{-\frac12}_{2,\infty}}
&\leq&C\|u_{n+m}-u_n\|_{B^{\frac12}_{2,1}}\|\partial_x\rho_{n+1}\|_{B^{-\frac12}_{2,\infty}}\nonumber\\
&\leq&C\|u_{n+m}-u_n\|_{B^{\frac12}_{2,1}}\|\rho_{n+1}\|_{B^{\frac12}_{2,1}}
\leq C\|u_{n+m}-u_n\|_{B^{\frac12}_{2,1}}.
\end{eqnarray}
For $R^1_{n,m}$, applying Lemma \ref{criproduct} to its decomposition in Step 3 of the proof of Theorem
\ref{LWP} and using the uniform bound (\ref{cribound}), we obtain
\begin{eqnarray}\label{criR1}
\|R^1_{n,m}\|_{B^{-\frac12}_{2,\infty}}
\leq C\Big(\|u_{n+m}-u_n\|_{B^{\frac12}_{2,1}}
+\|\rho_{n+m}-\rho_n\|_{B^{-\frac12}_{2,\infty}}\Big).
\end{eqnarray}
For the second equation of (\ref{diff}), the embedding $B^{\frac12}_{2,1}\hookrightarrow B^{\frac12}_{2,\infty}$
and the algebra property of $B^{\frac12}_{2,1}$ yield
\begin{eqnarray}\label{crid2}
\|(u_n-u_{n+m})\partial_xu_{n+1}\|_{B^{\frac12}_{2,\infty}}
\leq C\|(u_n-u_{n+m})\partial_xu_{n+1}\|_{B^{\frac12}_{2,1}}
\leq C\|u_{n+m}-u_n\|_{B^{\frac12}_{2,1}}.
\end{eqnarray}
It remains to estimate $R^2_{n,m}$. In light of
$\|\partial_xG\ast f\|_{B^{\frac12}_{2,\infty}}\leq C\|f\|_{B^{-\frac12}_{2,\infty}}$,
we only need to bound $F(u_{n+m},\rho_{n+m})-F(u_n,\rho_n)$ in $B^{-\frac12}_{2,\infty}$.
For the quartic term, we have
\begin{eqnarray}\label{crid4}
&&\|u_{n+m}^2\rho_{n+m}^2-u_n^2\rho_n^2\|_{B^{-\frac12}_{2,\infty}}\nonumber\\
&\leq&C\Big(\|u_{n+m}^2-u_n^2\|_{B^{\frac12}_{2,1}}\|\rho_{n+m}^2\|_{B^{-\frac12}_{2,\infty}}
+\|u_n^2\|_{B^{\frac12}_{2,1}}\|\rho_{n+m}^2-\rho_n^2\|_{B^{-\frac12}_{2,\infty}}\Big)\nonumber\\
&\leq&C\Big(\|u_{n+m}-u_n\|_{B^{\frac12}_{2,1}}+\|\rho_{n+m}-\rho_n\|_{B^{-\frac12}_{2,\infty}}\Big).
\end{eqnarray}
The remaining quadratic and cubic terms can be estimated in the same way. Consequently,
\begin{eqnarray}\label{criR2}
\|R^2_{n,m}\|_{B^{\frac12}_{2,\infty}}
\leq C\Big(\|u_{n+m}-u_n\|_{B^{\frac12}_{2,1}}
+\|\rho_{n+m}-\rho_n\|_{B^{-\frac12}_{2,\infty}}\Big).
\end{eqnarray}
Combining (\ref{crid1})--(\ref{criR2}) yields (\ref{cridiff}). By Lemma \ref{log}
and (\ref{cribound}), we infer that
\begin{eqnarray}\label{crilog}
\|u_{n+m}-u_n\|_{B^{\frac12}_{2,1}}
&\leq&C\|u_{n+m}-u_n\|_{B^{\frac12}_{2,\infty}}
\ln\Big(e+\frac{\|u_{n+m}-u_n\|_{B^{\frac32}_{2,1}}}{\|u_{n+m}-u_n\|_{B^{\frac12}_{2,\infty}}}\Big)\nonumber\\
&\leq&CA_{n,m}(t)\ln\Big(e+\frac{C}{A_{n,m}(t)}\Big)
\leq CA_{n,m}(t)\Big(1-\ln\frac{A_{n,m}(t)}{C}\Big),
\end{eqnarray}
where we used the fact that $x\mapsto x\ln(e+\frac cx)$ is
nondecreasing for any $c>0$. Plugging (\ref{crilog}) into (\ref{cridiff}) yields that
\begin{eqnarray}\label{criArec}
A_{n+1,m}(t)\leq C\Big(a_{n,m}+\int_0^tA_{n,m}(\tau)\Big(1-\ln\frac{A_{n,m}(\tau)}{C}\Big)d\tau\Big).
\end{eqnarray}
Let $A_n(t):=\sup_mA_{n,m}(t)$ and $a_n:=\sup_ma_{n,m}$. Enlarging $C$ if necessary, we may assume that
$A_n(t)<C$ on $[0,T]$. Since $x\mapsto x(1-\ln\frac{x}{C})$ is nondecreasing on
$[0,C)$, we deduce from (\ref{criArec}) that
\begin{eqnarray}\label{criArec2}
A_{n+1}(t)\leq C\Big(a_n+\int_0^tA_n(\tau)\Big(1-\ln\frac{A_n(\tau)}{C}\Big)d\tau\Big).
\end{eqnarray}
Note that $a_n\rightarrow0$ as $n\rightarrow\infty$, since $S_j\rho_0\rightarrow\rho_0$ in $B^{\frac12}_{2,1}$ and
$S_ju_0\rightarrow u_0$ in $B^{\frac32}_{2,1}$. Hence $\tilde A(t):=\limsup_{n\rightarrow\infty}A_n(t)$ satisfies
\begin{equation*}
\tilde A(t)\leq C\int_0^t\tilde A(\tau)\Big(1-\ln\frac{\tilde A(\tau)}{C}\Big)d\tau,
\end{equation*}
which together with Lemma \ref{Osgood}
implies $\tilde A(t)\equiv0$. Therefore, $\{(\rho_n,u_n)\}$ is a Cauchy sequence in
$L^\infty(0,T;B^{-\frac12}_{2,\infty})\times L^\infty(0,T;B^{\frac12}_{2,\infty})$ and converges to some limit
function $(\rho,u)$. The Fatou property and (\ref{cribound}) then ensure that
$(\rho,u)\in L^\infty(0,T;B^{\frac12}_{2,1})\times L^\infty(0,T;B^{\frac32}_{2,1})$. Passing to the limit in
(\ref{appro}), we deduce that $(\rho,u)$ is indeed a solution of (\ref{1.1}), and the transport equation theory
\cite{BCD} further yields that $(\rho,u)$ belongs to the class stated in Theorem \ref{LWPcri}.\\
\textbf{Step 4: Uniqueness and stability.}
Suppose that $(\rho_1,u_1)$ and $(\rho_2,u_2)$ are two solutions of (\ref{1.1}) and set
$A(t):=\|\rho_1(t)-\rho_2(t)\|_{B^{-\frac12}_{2,\infty}}+\|u_1(t)-u_2(t)\|_{B^{\frac12}_{2,\infty}}$.
We also set
\begin{equation*}
D(t):=\|\rho_1(t)-\rho_2(t)\|_{B^{-\frac12}_{2,\infty}}
+\|u_1(t)-u_2(t)\|_{B^{\frac12}_{2,1}}.
\end{equation*}
For the source terms in (\ref{uniq}), the same estimates as in (\ref{crid1})--(\ref{criR1})
give
\begin{eqnarray}\label{criuni1}
&&\|(u_2-u_1)\partial_x\rho_2\|_{B^{-\frac12}_{2,\infty}}\nonumber\\
&&\quad+\|\big(3\gamma_1u_1^2-(b-1)\big)\rho_1\partial_xu_1
-\big(3\gamma_1u_2^2-(b-1)\big)\rho_2\partial_xu_2\|_{B^{-\frac12}_{2,\infty}}
\leq CD(t).
\end{eqnarray}
Similarly, we deduce from (\ref{crid2})--(\ref{criR2}) that
\begin{eqnarray}\label{criuni2}
&&\|\alpha(u_2-u_1)\partial_xu_2\|_{B^{\frac12}_{2,\infty}}
+\|\partial_xG\ast\big(F(u_1,\rho_1)-F(u_2,\rho_2)\big)\|_{B^{\frac12}_{2,\infty}}
\leq CD(t).
\end{eqnarray}
Applying Lemma \ref{transport} to (\ref{uniq}) and using (\ref{criuni1})--(\ref{criuni2}), we obtain
\begin{equation*}
A(t)\leq C\Big(A(0)+\int_0^tD(\tau)d\tau\Big).
\end{equation*}
By Lemma \ref{log} and the uniform bounds of the two solutions, the same logarithmic argument as in
(\ref{crilog}) yields, with $C$ chosen large enough that $A(t)<C$,
\begin{equation*}
D(t)\leq CA(t)\Big(1-\ln\frac{A(t)}{C}\Big),
\end{equation*}
and hence
\begin{equation*}
A(t)\leq C\Big(A(0)+\int_0^tA(\tau)\Big(1-\ln\frac{A(\tau)}{C}\Big)d\tau\Big).
\end{equation*}
Applying Remark \ref{Osgoodrem}, we get
\begin{equation}\label{criuniq}
A(t)\leq C\,A(0)^{\exp\{-CT\}},\qquad \forall t\in[0,T).
\end{equation}
In particular, $A(0)=0$ implies $A(t)\equiv0$, which proves the uniqueness. Finally, with $\theta=\frac32-s'$, we have the interpolation
inequalities
\begin{eqnarray*}
\|\rho_1-\rho_2\|_{B^{s'-1}_{2,1}}
&\leq&C\|\rho_1-\rho_2\|_{B^{-\frac12}_{2,\infty}}^{\theta}
\|\rho_1-\rho_2\|_{B^{\frac12}_{2,1}}^{1-\theta},\\
\|u_1-u_2\|_{B^{s'}_{2,1}}
&\leq&C\|u_1-u_2\|_{B^{\frac12}_{2,\infty}}^{\theta}
\|u_1-u_2\|_{B^{\frac32}_{2,1}}^{1-\theta}.
\end{eqnarray*}
Combining these inequalities with (\ref{criuniq}) and the uniform bounds of the two solutions, we obtain
(\ref{holdercri}). This completes the proof of Theorem \ref{LWPcri}.
\end{proof}

\begin{remark3}\label{Sobolev}
$(i)$ Since $B^s_{2,2}=H^s$, taking $p=r=2$ in Theorem \ref{LWP} yields the local well-posedness of (\ref{1.1}) in
$C([0,T);H^{s-1})\times C([0,T);H^{s})$ with $s>\frac32$.\\
$(ii)$ For $b=2$ and $\tilde\beta=-\frac34$, the system (\ref{1.1}) reduces to the system (\ref{2CH}) modelling
equatorial shallow-water waves with the Coriolis effect \cite{HL}. Hence Theorem \ref{LWP} with $p=r=2$ recovers
the local well-posedness result of Hu and Liu \cite{HL} for the system (\ref{2CH}) in $H^{s-1}\times H^s$ with
$s>\frac32$. Moreover, Theorem \ref{LWPcri} also applies to the system (\ref{2CH}), and further yields its
local well-posedness in the critical Besov space $B^{\frac12}_{2,1}\times B^{\frac32}_{2,1}$.
\end{remark3}

\section{Blow-up criteria}\label{sec4}
\newtheorem{theorem4}{Theorem}[section]
\newtheorem{lemma4}{Lemma}[section]
\newtheorem {remark4}{Remark}[section]
\newtheorem {definition4}{Definition}[section]
\newtheorem{corollary4}{Corollary}[section]
\newtheorem{example4}{Example}[section]
\par
In this section, we first establish a general blow-up criterion for the strong solutions to (\ref{1.1}),
and then present a finite-time blow-up result for a special parameter choice. Our first result reads as follows.
\begin{theorem4}\label{blowup}
Let $(\rho_0,u_0)\in B^{\frac12}_{2,1}\times B^{\frac32}_{2,1}$ and let $T>0$ be the maximal existence time
of the corresponding solution $(\rho,u)$ to (\ref{1.1}). Then the solution blows up in finite time if and only if
\begin{eqnarray}\label{blowupcri}
\int_0^T\Big(\|u(t)\|_{L^\infty}^2+\|u_x(t)\|_{L^\infty}^2+\|\rho(t)\|_{L^\infty}^2
+\|u(t)\|_{L^\infty}^4+\|\rho(t)\|_{L^\infty}^4\Big)dt=\infty.
\end{eqnarray}
\end{theorem4}
\begin{proof}
We first derive an a priori estimate for the solution in the critical Besov space. Set
\begin{equation*}
N(t):=\|\rho(t)\|_{B^{\frac12}_{2,1}}+\|u(t)\|_{B^{\frac32}_{2,1}},
\end{equation*}
and
\begin{equation*}
Z(t):=1+\|u(t)\|_{L^\infty}^2+\|u_x(t)\|_{L^\infty}^2+\|\rho(t)\|_{L^\infty}^2
+\|u(t)\|_{L^\infty}^4+\|\rho(t)\|_{L^\infty}^4.
\end{equation*}
Applying $\Delta_j$ to the transport form (\ref{nonlocal}) yields
\begin{equation}\label{dyadic}
\left\{
 \begin{aligned}
&\partial_t\Delta_j\rho+u\partial_x\Delta_j\rho
=\Delta_j\Big(\big(3\gamma_1u^2-(b-1)\big)\rho u_x\Big)+R_j^1,\\
&\partial_t\Delta_ju+\alpha u\partial_x\Delta_ju
=-\Delta_j\partial_xG\ast F(u,\rho)+R_j^2,
\end{aligned}
\right.
\end{equation}
where
\begin{equation*}
R_j^1=[u\partial_x,\Delta_j]\rho,\qquad R_j^2=\alpha[u\partial_x,\Delta_j]u.
\end{equation*}
Multiplying the first equation of (\ref{dyadic}) by $\Delta_j\rho$ and the second one by $\Delta_ju$,
respectively, integrating over $\R$, and using integration by parts together with H\"older's inequality,
we obtain
\begin{eqnarray}\label{Ljrho}
\|\Delta_j\rho(t)\|_{L^2}
&\leq&\|\Delta_j\rho_0\|_{L^2}+C\int_0^t\Big(\|u_x\|_{L^\infty}\|\Delta_j\rho\|_{L^2}\nonumber\\
&&\quad+\|\Delta_j\big((3\gamma_1u^2-(b-1))\rho u_x\big)\|_{L^2}
+\|R_j^1\|_{L^2}\Big)d\tau
\end{eqnarray}
and
\begin{eqnarray}\label{Lju}
\|\Delta_ju(t)\|_{L^2}
\leq\|\Delta_ju_0\|_{L^2}+C\int_0^t\Big(\|u_x\|_{L^\infty}\|\Delta_ju\|_{L^2}+\|\Delta_j\partial_xG\ast F(u,\rho)\|_{L^2}+\|R_j^2\|_{L^2}\Big)d\tau.
\end{eqnarray}
We next multiply (\ref{Ljrho}) by $2^{\frac j2}$ and (\ref{Lju}) by $2^{\frac{3j}2}$, respectively,
sum over $j\geq-1$, and use Lemma \ref{commutator} to handle the commutator terms $R_j^1$ and $R_j^2$, obtaining
\begin{eqnarray}\label{RjrhoRju}
\big\|\big(2^{\frac j2}\|R_j^1\|_{L^2}\big)_j\big\|_{l^1}
\leq C\|u_x\|_{L^\infty}\|\rho\|_{B^{\frac12}_{2,1}},\qquad \big\|\big(2^{\frac{3j}2}\|R_j^2\|_{L^2}\big)_j\big\|_{l^1}
\leq C\|u_x\|_{L^\infty}\|u\|_{B^{\frac32}_{2,1}}.
\end{eqnarray}
Plugging (\ref{RjrhoRju}) into the above dyadic estimates, we arrive at
\begin{eqnarray}\label{rhoap}
\|\rho(t)\|_{B^{\frac12}_{2,1}}
\leq\|\rho_0\|_{B^{\frac12}_{2,1}}+C\int_0^t\Big(\|u_x\|_{L^\infty}\|\rho\|_{B^{\frac12}_{2,1}}
+\|\big(3\gamma_1u^2-(b-1)\big)\rho u_x\|_{B^{\frac12}_{2,1}}\Big)d\tau,
\end{eqnarray}
and
\begin{eqnarray}\label{uap}
\|u(t)\|_{B^{\frac32}_{2,1}}
\leq\|u_0\|_{B^{\frac32}_{2,1}}+C\int_0^t\Big(\|u_x\|_{L^\infty}\|u\|_{B^{\frac32}_{2,1}}
+\|F(u,\rho)\|_{B^{\frac12}_{2,1}}\Big)d\tau.
\end{eqnarray}
By Corollary \ref{algebra} and the embeddings
$B^{\frac12}_{2,1}\hookrightarrow L^\infty$ and
$B^{\frac32}_{2,1}\hookrightarrow B^{\frac12}_{2,1}$, all the nonlinear terms in (\ref{rhoap})--(\ref{uap})
can be controlled by $Z(t)N(t)$. For instance, the quartic term satisfies
\begin{eqnarray}\label{quarticap}
\|u^2\rho^2\|_{B^{\frac12}_{2,1}}
\leq C\Big(\|u\|_{L^\infty}^2\|\rho\|_{L^\infty}
+\|u\|_{L^\infty}\|\rho\|_{L^\infty}^2\Big)N(t)
\leq CZ(t)N(t).
\end{eqnarray}
The other terms can be estimated in the same way. Consequently,
\begin{eqnarray}\label{Nonlinest}
\|\big(3\gamma_1u^2-(b-1)\big)\rho u_x\|_{B^{\frac12}_{2,1}}
+\|F(u,\rho)\|_{B^{\frac12}_{2,1}}\leq CZ(t)N(t).
\end{eqnarray}
Combining (\ref{rhoap})--(\ref{Nonlinest}) yields
\begin{equation*}
N(t)\leq N(0)+C\int_0^tZ(\tau)N(\tau)d\tau,
\end{equation*}
and hence Gronwall's inequality gives
\begin{eqnarray}\label{Ngron}
N(t)\leq N(0)\exp\Big(C\int_0^tZ(\tau)d\tau\Big),\qquad \forall t\in[0,T).
\end{eqnarray}
Suppose that $T<\infty$ and the integral in (\ref{blowupcri}) is finite. Since $T<\infty$, this also gives
$\int_0^TZ(\tau)d\tau<\infty$, and (\ref{Ngron}) implies
\begin{equation*}
\limsup_{t\to T}\Big(\|\rho(t)\|_{B^{\frac12}_{2,1}}+\|u(t)\|_{B^{\frac32}_{2,1}}\Big)<\infty.
\end{equation*}
The local existence theory and uniqueness then allow the solution to be continued beyond $T$, contradicting
the maximality of $T$. Therefore, the integral in (\ref{blowupcri}) must diverge when $T<\infty$.

Conversely, the boundedness of $N(t)$ for any strong solution on $[0,T]$
implies that the integral in (\ref{blowupcri}) is finite. Hence, if this
integral diverges, the solution blows up in finite time. This completes the proof of Theorem \ref{blowup}.
\end{proof}

Motivated by \cite{CE1}, we next consider the special case $b=3$ and $\tilde\beta=\frac12$, for which
\begin{equation*}
\alpha=\frac74,\qquad \gamma_1=18,\qquad \gamma_2=8,\qquad
\gamma_3=-4\sqrt2,\qquad \gamma_4=0.
\end{equation*}
In this case, the condition $\gamma_4=0$ ensures that the symmetry
$(\rho,u)\mapsto(\rho(t,-x),-u(t,-x))$ is preserved, which yields the following finite-time blow-up result.
\begin{theorem4}\label{symmetryblowup}
Let $b=3$, $\tilde\beta=\frac12$ and $(\rho_0,u_0)\in H^{s-1}\times H^s$ with $s>\frac32$.
Let $T$ be the maximal existence time of the corresponding solution $(\rho,u)$ to (\ref{1.1}).
Assume that $u_0$ is odd, $\rho_0$ is even, $u_0'(0)<0$ and $\rho_0(0)=0$. Then
\begin{equation}\label{T0}
T\leq T_0:=-\frac{1}{\alpha u_0'(0)}=-\frac{4}{7u_0'(0)}.
\end{equation}
Moreover, for every $t\in[0,T)$,
\begin{equation}\label{slopebound}
u_x(t,0)\leq\frac{u_0'(0)}{1+\alpha u_0'(0)t}.
\end{equation}
In particular, the solution blows up in finite time.
\end{theorem4}
\begin{proof}
By a density argument, we may assume $s=3$. Since $u_0$ is odd, $\rho_0$ is even and $\gamma_4=0$,
the structure of (\ref{1.1}) and uniqueness imply that $u(t,\cdot)$ is odd and $\rho(t,\cdot)$ is even
for all $t\in[0,T)$. Consequently,
\begin{equation}\label{originsymmetry}
u(t,0)=0,\qquad u_{xx}(t,0)=0,\qquad \rho_x(t,0)=0.
\end{equation}

Set $M(t):=u_x(t,0)$. Using the first equation of (\ref{1.1}) and (\ref{originsymmetry}), we obtain
\begin{equation*}
\frac{d}{dt}\rho(t,0)=\Big(-2\rho u_x+54u^2u_x\rho\Big)(t,0)=-2\rho(t,0)M(t),
\end{equation*}
and hence
\begin{equation}\label{rhoorigin}
\rho(t,0)=\rho_0(0)\exp\Big(-2\int_0^tM(s)ds\Big)=0,\qquad \forall t\in[0,T).
\end{equation}
On the other hand, differentiating the second equation of (\ref{nonlocal}) with respect to $x$ and using
$\partial_x^2G\ast F=G\ast F-F$ yield
\begin{equation}\label{slopeeq}
u_{xt}+\alpha u_x^2+\alpha uu_{xx}=F(u,\rho)-G\ast F(u,\rho).
\end{equation}
For the present parameters, $F(u,\rho)$ can be rewritten as
\begin{equation}\label{Fsquare}
F(u,\rho)=\frac{13}{8}u^2+\frac12\rho^2
+8u^2\Big(\rho-\frac{\sqrt2}{4}\Big)^2\geq0.
\end{equation}
Evaluating (\ref{slopeeq}) at $x=0$ and using (\ref{originsymmetry}), (\ref{rhoorigin}) and (\ref{Fsquare}), we deduce that
\begin{equation}\label{Riccati}
\frac{dM}{dt}=-\alpha M^2-G\ast F(u,\rho)(t,0)\leq-\alpha M^2,
\end{equation}
because $G(x)=\frac12e^{-|x|}>0$ and $F(u,\rho)\geq0$.
Since $M(0)=u_0'(0)<0$, the differential inequality (\ref{Riccati}) implies $M(t)\leq M(0)<0$,
and a comparison with the Riccati equation gives
\begin{equation}\label{Mbound}
M(t)\leq\frac{M(0)}{1+\alpha M(0)t}
=\frac{u_0'(0)}{1+\alpha u_0'(0)t},\qquad 0\leq t<\min\{T,T_0\}.
\end{equation}
The right-hand side tends to $-\infty$ as $t\uparrow T_0$. If $T>T_0$, then $M(t)$ would remain bounded
on $[0,T_0]$, contrary to (\ref{Mbound}). Therefore $T\leq T_0$, and (\ref{slopebound}) holds on the whole
interval $[0,T)$. This completes the proof.
\end{proof}

\begin{remark4}\label{symmetryremark}
A typical choice of initial data satisfying the assumptions of Theorem \ref{symmetryblowup} is
\begin{equation*}
u_0(x)=-xe^{-x^2},\qquad \rho_0(x)=x^2e^{-x^2}.
\end{equation*}
Indeed, $u_0$ is odd, $\rho_0$ is even, $u_0'(0)=-1<0$ and $\rho_0(0)=0$.
Thus the corresponding solution blows up in finite time with $T\leq T_0=\frac47$.
\end{remark4}

\appendix
\stepcounter{section}
\section*{Appendix}\label{sec5}
\newtheorem{lemma5}{Lemma}[section]
\newtheorem {remark5}{Remark}[section]
\newtheorem{corollary5}{Corollary}[section]
\newtheorem{proposition5}{Proposition}[section]
\par
In this appendix, we recall some facts on the Littlewood-Paley decomposition, the nonhomogeneous Besov spaces
and the transport equation theory, which have been used in the previous sections. For more details, one can
refer to \cite{BCD,D}.

We first recall the Littlewood-Paley decomposition. Let $\mathcal{C}$ be the annulus
$\{\xi\in\R:\frac34\leq|\xi|\leq\frac83\}$. There exist radial functions $\chi$ and $\varphi$, valued in the
interval $[0,1]$, belonging respectively to $\mathcal{D}(B(0,\frac43))$ and $\mathcal{D}(\mathcal{C})$, such that
\begin{equation*}
\forall\xi\in\R,\quad \chi(\xi)+\sum_{j\geq0}\varphi(2^{-j}\xi)=1,
\end{equation*}
\begin{equation*}
|j-j'|\geq2\Rightarrow\mbox{Supp}\ \varphi(2^{-j}\cdot)\cap\mbox{Supp}\ \varphi(2^{-j'}\cdot)=\emptyset,
\quad j\geq1\Rightarrow\mbox{Supp}\ \chi\cap\mbox{Supp}\ \varphi(2^{-j}\cdot)=\emptyset.
\end{equation*}
Denote by $\mathcal{F}$ the Fourier transform and $\mathcal{F}^{-1}$ its inverse, and write
$h=\mathcal{F}^{-1}\varphi$ and $\tilde{h}=\mathcal{F}^{-1}\chi$. The nonhomogeneous dyadic blocks $\Delta_j$
are defined by
\begin{eqnarray*}
&&\Delta_ju=0,\quad if\ j\leq-2,\quad\quad
\Delta_{-1}u=\chi(D)u=\int_\R\tilde{h}(y)u(x-y)dy,\\
&&\Delta_ju=\varphi(2^{-j}D)u=2^{j}\int_\R h(2^jy)u(x-y)dy,\quad if\ j\geq0,
\end{eqnarray*}
and the low-frequency cut-off operators $S_j$ are defined by
\begin{equation*}
S_ju=\sum_{j'\leq j-1}\Delta_{j'}u.
\end{equation*}
The nonhomogeneous Besov spaces are defined by
\begin{equation*}
B^s_{p,r}(\R):=\Big\{u\in\mathcal{S}'(\R):\|u\|_{B^s_{p,r}(\R)}:=
\Big(\sum_{j\geq-1}2^{rjs}\|\Delta_ju\|^r_{L^p(\R)}\Big)^{\frac1r}<\infty\Big\},
\end{equation*}
with the usual modification when $r=\infty$.

Next, we introduce some useful lemmas and propositions about Besov spaces.

\begin{proposition5}[\cite{BCD}]\label{embedding}
Let $1\leq p_1\leq p_2\leq\infty$ and $1\leq r_1\leq r_2\leq\infty$, and let $s$ be a real number. Then we have
\begin{equation*}
B^s_{p_1,r_1}(\R)\hookrightarrow B^{s-(\frac{1}{p_1}-\frac{1}{p_2})}_{p_2,r_2}(\R).
\end{equation*}
If $s>\frac1p$ or $s=\frac1p$, $r=1$, we then have $B^s_{p,r}(\R)\hookrightarrow L^\infty(\R)$.
Moreover, for any $s'<s''$ and $\theta\in(0,1)$, the following interpolation inequality holds:
\begin{equation*}
\|f\|_{B^{\theta s'+(1-\theta)s''}_{p,r}(\R)}
\leq C\|f\|_{B^{s'}_{p,r}(\R)}^{\theta}\|f\|_{B^{s''}_{p,r}(\R)}^{1-\theta}.
\end{equation*}
\end{proposition5}

\begin{corollary5}[\cite{BCD}]\label{algebra}
$(i)$ For any positive real number $s$ and any $(p,r)$ in $[1,\infty]^2$, the space
$L^\infty(\R)\cap B^s_{p,r}(\R)$ is an algebra, and a constant $C$ exists such that
\begin{equation*}
\|uv\|_{B^s_{p,r}(\R)}\leq C\big(\|u\|_{L^\infty(\R)}\|v\|_{B^s_{p,r}(\R)}+\|u\|_{B^s_{p,r}(\R)}\|v\|_{L^\infty(\R)}\big).
\end{equation*}
$(ii)$ If $s>\frac1p$ or $s=\frac1p$, $r=1$, then $B^s_{p,r}(\R)$ is an algebra, i.e.,
\begin{equation*}
\|uv\|_{B^s_{p,r}(\R)}\leq C\|u\|_{B^s_{p,r}(\R)}\|v\|_{B^s_{p,r}(\R)}.
\end{equation*}
\end{corollary5}

\begin{lemma5}[Morse-type estimate, \cite{BCD,D}]\label{Morse}
Let $\sigma>\max\{\frac1p,\frac12\}$ and $(p,r)\in[1,\infty]^2$. For any $a\in B^{\sigma-1}_{p,r}(\R)$ and
$b\in B^{\sigma}_{p,r}(\R)$, there exists a constant $C$ such that
\begin{equation*}
\|ab\|_{B^{\sigma-1}_{p,r}(\R)}\leq C\|a\|_{B^{\sigma-1}_{p,r}(\R)}\|b\|_{B^{\sigma}_{p,r}(\R)}.
\end{equation*}
\end{lemma5}

\begin{lemma5}[Commutator estimate, \cite{BCD}]\label{commutator}
Let $\sigma>0$, $1\leq r\leq\infty$ and $1\leq p\leq p_1\leq\infty$. Let $v$ be a vector field over $\R$.
Define $R_j=[v\partial_x,\Delta_j]f$. There exists a constant $C$ such that
\begin{equation*}
\big\|\big(2^{j\sigma}\|R_j\|_{L^p(\R)}\big)_j\big\|_{l^r}
\leq C\Big(\|\partial_xv\|_{L^\infty(\R)}\|f\|_{B^\sigma_{p,r}(\R)}
+\|\partial_xf\|_{L^{p_2}(\R)}\|\partial_xv\|_{B^{\sigma-1}_{p_1,r}(\R)}\Big),
\end{equation*}
where $\frac{1}{p_2}=\frac1p-\frac{1}{p_1}$. Further, if $\sigma<1$, then
\begin{equation*}
\big\|\big(2^{j\sigma}\|R_j\|_{L^p(\R)}\big)_j\big\|_{l^r}
\leq C\|\partial_xv\|_{L^\infty(\R)}\|f\|_{B^\sigma_{p,r}(\R)}.
\end{equation*}
\end{lemma5}

The following lemma is crucial to the study of well-posedness in the critical space
$B^{\frac12}_{2,1}(\R)\times B^{\frac32}_{2,1}(\R)$.

\begin{lemma5}[\cite{LY2}]\label{criproduct}
For any $a\in B^{-\frac12}_{2,\infty}(\R)$ and $b\in B^{\frac12}_{2,1}(\R)$, there exists a constant $C$ such that
\begin{equation*}
\|ab\|_{B^{-\frac12}_{2,\infty}(\R)}\leq C\|a\|_{B^{-\frac12}_{2,\infty}(\R)}\|b\|_{B^{\frac12}_{2,1}(\R)}.
\end{equation*}
\end{lemma5}

\begin{lemma5}[Logarithmic interpolation inequality, \cite{D}]\label{log}
For any $f\in B^{\frac32}_{2,1}(\R)$, there exists a constant $C$ such that
\begin{equation*}
\|f\|_{B^{\frac12}_{2,1}(\R)}
\leq C\|f\|_{B^{\frac12}_{2,\infty}(\R)}
\ln\Big(e+\frac{\|f\|_{B^{\frac32}_{2,1}(\R)}}{\|f\|_{B^{\frac12}_{2,\infty}(\R)}}\Big).
\end{equation*}
\end{lemma5}

\begin{remark5}\label{Fatou}
Let $s\in\R$, $1\leq p,r\leq\infty$. The following properties hold true (see \cite{BCD}):\\
$(i)$ $B^s_{p,r}(\R)$ is a Banach space and is continuously embedded into $\mathcal{S}'(\R)$, where
$\mathcal{S}'(\R)$ is the dual space of the Schwartz space $\mathcal{S}(\R)$.\\
$(ii)$ If $p,r<\infty$, then $\mathcal{S}(\R)$ is dense in $B^s_{p,r}(\R)$.\\
$(iii)$ (Fatou property) If $(u_n)$ is a bounded sequence of $B^s_{p,r}(\R)$, then an element
$u\in B^s_{p,r}(\R)$ and a subsequence $u_{n_k}$ exist such that
$\lim_{k\rightarrow\infty}u_{n_k}=u$ in $\mathcal{S}'(\R)$ and
$\|u\|_{B^s_{p,r}(\R)}\leq C\liminf_{k\rightarrow\infty}\|u_{n_k}\|_{B^s_{p,r}(\R)}$.\\
$(iv)$ $B^s_{2,2}(\R)=H^s(\R)$.
\end{remark5}

The following Osgood lemma appears as a substitution for Gronwall's lemma.

\begin{lemma5}[Osgood's lemma, \cite{BCD}]\label{Osgood}
Let $\rho\geq0$ be a measurable function, $\gamma>0$ be a locally integrable function and $\mu$ be a continuous
and increasing function. Assume that, for some nonnegative real number $c$, the function $\rho$ satisfies
\begin{equation*}
\rho(t)\leq c+\int_{t_0}^t\gamma(t')\mu(\rho(t'))dt'.
\end{equation*}
If $c>0$, then $-\mathcal{M}(\rho(t))+\mathcal{M}(c)\leq\int_{t_0}^t\gamma(t')dt'$ with
$\mathcal{M}(x)=\int_x^1\frac{dr}{\mu(r)}$.\\
If $c=0$ and $\mu$ satisfies the condition $\int_0^1\frac{dr}{\mu(r)}=+\infty$, then the function $\rho=0$.
\end{lemma5}

\begin{remark5}\label{Osgoodrem}
In this paper, we set $\mu(r)=r(1-\ln\frac{r}{C})$, which satisfies the condition $\int_0^1\frac{dr}{\mu(r)}=+\infty$.
A simple calculation shows that $\mathcal{M}(x)=\ln\big(1-\ln\frac{x}{C}\big)-\ln\big(1+\ln C\big)$-type function,
and then we deduce from Lemma \ref{Osgood} that
\begin{equation*}
\rho(t)\leq c^{\exp\big(-\int_{t_0}^t\gamma(t')dt'\big)}C^{1-\exp\big(-\int_{t_0}^t\gamma(t')dt'\big)},
\quad if\ c>0,
\end{equation*}
and $\rho(t)=0$ if $c=0$.
\end{remark5}

Finally, we introduce the a priori estimates for the following transport equation
\begin{equation}\label{transeq}
\left\{
 \begin{aligned}
 &f_t+v\partial_xf=g,\\
&f|_{t=0}=f_0.
\end{aligned}
\right.
\end{equation}

\begin{lemma5}[A priori estimates in Besov spaces, \cite{BCD}]\label{transport}
Let $1\leq p\leq p_1\leq\infty$, $1\leq r\leq\infty$, $s\geq-\min\{\frac{1}{p_1},1-\frac1p\}$.
For the solution $f\in L^\infty(0,T;B^s_{p,r}(\R))$ of (\ref{transeq}) with the velocity
$\partial_xv\in L^1(0,T;B^{s-1}_{p,r}(\R)\cap L^\infty(\R))$-type regularity, the initial data
$f_0\in B^s_{p,r}(\R)$ and $g\in L^1(0,T;B^s_{p,r}(\R))$, we have
\begin{eqnarray*}
&&\|f(t)\|_{B^s_{p,r}(\R)}\leq\|f_0\|_{B^s_{p,r}(\R)}
+\int_0^t\Big(\|g(t')\|_{B^s_{p,r}(\R)}+CV'_{p_1}(t')\|f(t')\|_{B^s_{p,r}(\R)}\Big)dt',\\
&&\|f\|_{L^\infty_t(B^s_{p,r}(\R))}\leq\Big(\|f_0\|_{B^s_{p,r}(\R)}
+\int_0^te^{-CV_{p_1}(t')}\|g(t')\|_{B^s_{p,r}(\R)}dt'\Big)e^{CV_{p_1}(t)},
\end{eqnarray*}
where $V_{p_1}(t)=\int_0^t\|\partial_xv\|_{B^{\frac{1}{p_1}}_{p_1,\infty}(\R)\cap L^\infty(\R)}dt'$ if
$s<1+\frac{1}{p_1}$, $V_{p_1}(t)=\int_0^t\|\partial_xv\|_{B^{s-1}_{p_1,r}(\R)}dt'$ if $s>1+\frac{1}{p_1}$ or
$s=1+\frac{1}{p_1}$, $r=1$, and $C$ is a constant depending only on $s,p,p_1$ and $r$.
\end{lemma5}

In the case $s=1+\frac1p$, we have the following a priori estimates.

\begin{lemma5}[\cite{LY2}]\label{transportcri}
Let $1\leq p\leq p_1\leq\infty$ and $1\leq r\leq\infty$. For the solution
$f\in L^\infty(0,T;B^{1+\frac1p}_{p,r}(\R))$ of (\ref{transeq}) with the velocity
$\partial_xv\in L^1(0,T;B^{\frac1p}_{p_1,r}(\R)\cap L^\infty(\R))$, the initial data
$f_0\in B^{1+\frac1p}_{p,r}(\R)$ and $g\in L^1(0,T;B^{1+\frac1p}_{p,r}(\R))$, we have
\begin{eqnarray*}
&&\|f(t)\|_{B^{1+\frac1p}_{p,r}(\R)}\leq\|f_0\|_{B^{1+\frac1p}_{p,r}(\R)}
+\int_0^t\Big(\|g(t')\|_{B^{1+\frac1p}_{p,r}(\R)}+CV'_{p_1}(t')\|f(t')\|_{B^{1+\frac1p}_{p,r}(\R)}\Big)dt',\\
&&\|f\|_{L^\infty_t(B^{1+\frac1p}_{p,r}(\R))}\leq\Big(\|f_0\|_{B^{1+\frac1p}_{p,r}(\R)}
+\int_0^te^{-CV_{p_1}(t')}\|g(t')\|_{B^{1+\frac1p}_{p,r}(\R)}dt'\Big)e^{CV_{p_1}(t)},
\end{eqnarray*}
where $V_{p_1}(t)=\int_0^t\|\partial_xv\|_{B^{\frac1p}_{p_1,\infty}(\R)\cap L^\infty(\R)}dt'$ and $C$ is a
constant depending only on $p,p_1$ and $r$.
\end{lemma5}

\noindent\textbf{Acknowledgments.}
The authors thank the anonymous referee for helpful suggestions and comments.
This work was supported by the National Natural Science Foundation of China under Grant 12001528.

\end{document}